\documentclass[a4paper,11pt,twoside]{article}
\usepackage{latexsym,amsmath,amssymb,amsfonts,amsbsy,fancyhdr}
\usepackage{theorem}
\usepackage{amscd}
\usepackage{enumerate}
\usepackage{mathrsfs}
\usepackage{dsfont}
\usepackage[active]{srcltx}
\usepackage[colorlinks=false, pdfstartview=FitV, linkcolor=red, citecolor=green, urlcolor=red]{hyperref}
\usepackage{nicefrac}
\usepackage{marvosym}
\usepackage{longtable}
\usepackage{subfigure}
\usepackage{color}

\newcommand{\R}{\mathbb R}

\newcommand{\C}{\mathbb C}
\newcommand{\N}{\mathbb N}

\newcommand{\E}{\mathbb E}
\renewcommand{\S}{\mathcal S}
\newcommand{\F}{\mathcal F}
\renewcommand{\L}{\mathscr L}

\newcommand{\dd}{{\rm{d}}}
\newcommand{\qg}{\mathcal{Q}_g}
\newcommand{\dt}{\frac{\partial}{\partial t}}
\newcommand{\du}{\frac{\partial}{\partial u}}
\newcommand{\ds}{\frac{\partial}{\partial s}}
\newcommand{\dr}{\frac{\partial}{\partial r}}
\newcommand{\ddt}{\frac{\dd }{\dd t}}

\newcommand{\we}{\exp^\diamond}

\newtheorem{theorem}{Theorem}

\newtheorem{proposition}[theorem]{Proposition}
\newtheorem{lemma}[theorem]{Lemma}
\theorembodyfont{\rm}
\newtheorem{definition}[theorem]{Definition}
\newtheorem{example}[theorem]{Example}

\newtheorem{remark}[theorem]{Remark}
\newenvironment{rem}{\begin{remark}}{\hfill$\lozenge$\end{remark}}
\newtheorem{proof}{Proof}

\newenvironment{pf}{\begin{proof}}{\hfill$\square$\end{proof}}

\allowdisplaybreaks

\DeclareMathAlphabet{\mathcal}{OMS}{cmsy}{m}{n}

\begin{document}
\title{A generalised It\=o formula for L\'evy-driven Volterra processes\footnote{Author manuscript, accepted for publication in {\it Stochastic Processes and their Applications},
\newline
DOI: 10.1016/j.spa.2015.02.009}}
\author{Christian Bender, Robert Knobloch and Philip Oberacker
\thanks{Department of Mathematics, Saarland University, P.O. Box 151150, 66041
Saarbr\"ucken, Germany
\newline
e-mail: bender@math.uni-sb.de, knobloch@math.uni-sb.de,
oberacker@math.uni-sb.de}}
\maketitle

\begin{abstract}
We derive a generalised It\=o formula for stochastic processes which are constructed by a convolution of a deterministic kernel with a centred L\'evy process.
This formula has a unifying character in the sense that it contains the classical It\=o formula for L\'evy processes as well as recent change-of-variable formulas for Gaussian
processes such as fractional Brownian motion as special cases. Our result also covers  fractional L\'evy processes
(with Mandelbrot-Van Ness kernel) and a wide class of  related processes for which
such a generalised It\=o formula has  not yet been available in the literature.
\end{abstract}

\noindent {\bf AMS 2010 Mathematics Subject Classification:} 60G22, 60H05, 60H07

\noindent {\bf Keywords:} Fractional L\'evy process, It\=o formula, Skorokhod integral, Stochastic convolution, \newline $S$-transform

\section{Introduction}

In this paper we consider processes that are constructed by a convolution of a deterministic Volterra kernel $f$ and a centred two-sided L\'evy process $L$, i.e.
$$
M(t)=\int_{-\infty}^t f(t,s) \ L(\dd s), \quad t\geq 0,
$$
which we refer to as L\'evy-driven Volterra processes. If $L$ has no jump part, this construction contains Volterra Gaussian processes such as fractional Brownian motion, for which a Skorokhod type stochastic calculus beyond the semimartingale case
was developed in recent years, see e.g. \cite{AMN01}, \cite{B_sta} and \cite{HO03} for different approaches.

However, in many applications such as financial modelling, see e.g. \cite{BBV_mesp} or \cite{KM_gfl}, signal processing \cite{Ual14} or network traffic \cite{WT05}, the restriction to Gaussian distributions cannot always be justified. As the covariance
structure of
L\'evy-driven Volterra processes is determined by the choice of the kernel, these processes have the same long memory or short memory properties as the corresponding Gaussian ones, but allow for more
flexibility in the choice of the distribution. 

Fractional L\'evy processes, where $f$ is the Mandelbrot-Van Ness kernel of a fractional Brownian motion, are among the best studied 
L\'evy-driven Volterra processes \cite{EW_uaflp,Mar06,Neu14,TM11}. It is known that fractional L\'evy processes and more general classes of L\'evy-driven Volterra processes may fail to be semimartingales \cite{BP_ldma, BLS12}, and so
the classical It\=o calculus does not apply to these processes. A first step towards a Skorokhod type integration theory for L\'evy-driven Volterra processes was done in \cite{BM_sccl}, but only
in the case of a pure jump L\'evy process and under very restrictive 
assumptions on the kernel and strong moment conditions on the L\'evy measure. In the present paper we first provide a general framework for a Skorokhod integration with respect 
to L\'evy-driven Volterra processes, when the driving L\'evy process is square-integrable. We follow the $S$-transform approach of \cite{B_sta,BM_sccl,Jan97} which avoids some technicalities of Malliavin calculus by restricting
the duality relation between Malliavin derivative and Skorokhod integral to an appropriate set of stochastic exponentials as test functionals. The Hitsuda-Skorokhod integral with respect to $M$ 
introduced in Definition~\ref{M_diamond} below extends the classical It\=o integral for L\'evy processes  and also covers divergence type integrals for Gaussian processes. We note that our 
Hitsuda-Skorokhod integral with respect to L\'evy-driven Volterra processes is closely related to the integral, which was developed in the recent papers  \cite{B-NBPV_sivm,B-NBS_sivm,BS_it}
for volatility modulated L\'evy-driven Volterra processes in the cases that either $L$ has no jumps or is of pure jump type. Indeed, the main difference is that in the latter papers a Malliavin trace term is added, which under suitable assumptions 
changes from the Hitsuda-Skorokhod integral to the backward integral, if $L$ has no jump part.

As a main result we derive -- under appropriate assumptions on the kernel and the driving L\'evy process -- the following generalised It\=o formula:
\begin{equation}
 \begin{aligned}\label{ito_formula_L2}
G(M(T)) &= G(0) + \frac{\sigma^2}{2} \int_{0}^T G''(M(t)) \left( \ddt
\int_{-\infty}^t f(t,s)^2  \dd s \right) \dd t\\
&\quad+ \sum_{0<t\leq T} \left[ G(M(t))- G(M(t-))- G'(M(t-)) \Delta M(t)
\right] 
\\
&\quad+  \int_{-\infty}^T \int_{\R_0} \int_{0 \lor s}^T \left( G'\left(M(t) +
xf(t,s)\right) -G'(M(t)) \right) x \dt f(t,s) \ \dd t \ N^\diamond(\dd x, \dd s)  \\
&\quad+ \int_{0}^T G'(M(t)) \ M^\diamond (\dd t).
 \end{aligned}
\end{equation}
Here $N$ denotes the jump measure and $\sigma$ the standard deviation of the Gaussian part of the L\'evy process $L$. Moreover, $G$ is twice continuously 
differentiable and satisfies some growth conditions. The diamonds indicate that integration is understood in the Skorokhod sense. Let us stress that this formula
is unifying in the following sense:
\begin{itemize}
 \item If we take $f(t,s)=\mathds1_{[0,t]}(s)$, then $M$ coincides with the driving L\'evy process $L$, the triple integral on the right-hand side vanishes, and 
we recover the classical It\=o formula for L\'evy processes.
\item If $L$ has no jump part, then $M$ is a Gaussian process. In this case
the triple integral and the sum over the jumps of $M$ vanish. Moreover, the expression $\sigma^2 \int_{-\infty}^t f(t,s)^2  \dd s$ then
equals the variance of $M$. Hence, we end up with the It\=o formula for Gaussian processes in the Wick-It\=o sense, which has been proved in the literature by various techniques, see e.g. \cite{AMN01,Ben03,Ben14,BOSW04,NT06}.
\item If $L$ has no Gaussian part, the first integral on the right-hand side vanishes, and such a formula was proved by \cite{BM_sccl} when the kernel $f$ has compact support and sufficient smoothness
and the L\'evy process has finite moments of all orders. In this case we here considerably weaken the integrability assumptions on the L\'evy process and extend the class of kernels for which 
the generalised It\=o formula holds. In particular, our result now covers the fractional L\'evy case which was excluded in \cite{BM_sccl}.
\end{itemize}
In the presence of a Gaussian part we also show that the growth conditions on $G$ can be relaxed, which is due to the rapid decrease of the characteristic function of $M$ 
(under a change of measure) in this case.

The paper is organised as follows: In Section 2 we describe the set-up and in particular introduce the assumptions on the kernel $f$. Section 3 is devoted to the definition and discussion
of the Hitsuda-Skorokhod integral with respect to the L\'evy-driven Volterra process $M$ and with respect to the jump measure. The proof of the above generalised It\=o formula is given in Section 4.

 \section{Set-up}\label{s.setup}
This section aims at introducing the Volterra type processes the present paper is concerned with. In particular, here we introduce the notation regarding the underlying L\'evy process. For more information on the theory of L\'evy processes we refer e.g. to
\cite{Ber96} and \cite{Kyp06} as well as to \cite{App09}, where stochastic analysis with respect to L\'evy
processes is treated.

Throughout this paper we fix some $T \ge 0$.

Let $(L_1(t))_{t \geq 0}$ and $(L_2(t))_{t \geq 0}$ be two independent L\'evy processes on a complete probability space $(\Omega,\mathscr F,\mathbb P)$ with characteristic triplet   $(\gamma, \sigma, \nu)$, where $\sigma\ge0$,  $\nu$ is a L\'evy measure on $\R_0$, that satisfies 

\begin{equation}\label{e.L-k}
\int_{\R_0} x^2 \ \nu(\dd x) <\infty,
\end{equation}
and
\begin{equation}\label{drift_gamma}
  \gamma=-\int_{\R\setminus[-1,1]}x \ \nu(\dd x).
\end{equation}

For $l=1,2$ note that (\ref{e.L-k}) and (\ref{drift_gamma}) are equivalent to  $L_l(t)\in\mathscr L^2(\mathbb P)$ and  $\E(L_l(1))=0$. Moreover, the processes ${L_{l}}(t)$ can be represented as

\[
L_{l}(t) =  \sigma W_l(t)+ \int_0^t \int_{\R_0} x \ \tilde{N_l}(\dd x, \dd s),
\]

where $ W_l$ is a standard Brownian motion and $ \tilde{N_l}(\dd x, \dd s) =  N_l(\dd x, \dd s)- \nu(\dd x) \dd s$ is the compensated jump measure of the L\'evy process $L_l$.

We define a two-sided L\'evy process $L:=(L(t))_{t \in \R}$  by

\[
 L(t) := \begin{cases} L_1(t), &  t \geq 0\\-L_2(-(t-)), & t<0
 \end{cases}
\] 

and the two-sided Brownian motion $W:=(W(t))_{t \in \R}$ is defined analogously.

Note that for $a,b \in \R$ with $a \le b$ and any Borel set $\mathcal A \subset \R_0$ the jump measure $N(\dd x, \dd s)$ of the two-sided process $L$ fulfils 

\[
 \begin{aligned}
  N(\mathcal A, [a,b]) &= \# \{s \in [a,b]: \ \Delta L(s) \in  \mathcal A\} \\
  &= \# \{s \in [a,b] \cap [0,\infty): \ \Delta L_1(s) \in  \mathcal A\} \ + \ \# \{s \in [a,b] \cap (-\infty,0): \ \Delta L_2(-s) \in  \mathcal A\}
 \end{aligned}
\]

and

\[
 \E\left(N(\mathcal A,[a,b])\right) = (b-a) \, \nu(\mathcal A).
\]

The compensated jump measure of the two-sided process $L$ is defined as $\tilde N(\dd x, \dd s) := N(\dd x, \dd s) - \nu(\dd x) \ \dd s$.
Furthermore, we assume that $\mathscr F$ is the completion of the $\sigma$-algebra generated by $L$.

We now introduce the main object under consideration in the present paper, which is the following class of stochastic integrals with respect to $L$.

\begin{definition}\label{d.2}
For any  function $f:\R^2 \to \R$ such that $f(t,\cdot) \in \L^2(\dd s)$ for every $t \in \R$ we define a stochastic process $M:=(M(t))_{t\ge 0}$ by
\[
M(t)=\int_{-\infty}^t f(t,s) \ L(\dd s)
\]
for every $t\ge 0$. The process $M$ is referred to as {\it L\'evy-driven Volterra process}.
\end{definition}

Throughout this paper we shall work with a class $\mathcal K$ of integration kernels which is defined as follows:
\begin{definition}\label{d.1}
We denote by $\mathcal K$ the class of measurable functions $f:\R^2\to\R$ with $\text{supp}f \subset [\tau, \infty)^2$ for some $\tau \in [-\infty, 0]$ such that
\begin{enumerate}[(i)]
\item\label{volterra} $\forall\,s>t\ge0:\quad f(t,s)=0$,
\item\label{f(0)} $f(0,\cdot)=0\quad$ Lebesgue-a.e.,

    \item\label{continuity_of_f} the function $f$ is continuous on the set $\{(t,s) \in \R^2: \  \tau \le s \le t \}$,
    
     \item\label{no_zero_set} for all $0< t  \le T$, $\{s\in\R:f(t,s)\ne0\}$ is not a Lebesgue null set,

    \item\label{behaviour_kernel}
    for all $s \in \R$ the mapping $t \mapsto  f(t,s)$ is continuously differentiable on the set $(s,\infty)$ and there exist some $C_0>0$ and $\beta, \gamma \in (0,1)$ with $\beta + \gamma <1$ such that
    \begin{equation}\label{e.decay_of_dtf}
      \left|\dt f(t,s) \right| \le C_0|s|^{-\beta} |t-s|^{-\gamma}
    \end{equation}
    for all $t>s>-\infty$. Furthermore, for any $t \in [0,T]$ there exists some $\epsilon >0$ such that 
    
    \begin{equation}\label{e.decay_of_f}
   \sup_{s \in (-\infty,-1]} \left( \sup_{r \in [0 \lor(t-\epsilon),\, (t+\epsilon) \land T]}|f(r,s)| |s|^{\theta} \right) < \infty
    \end{equation}
    
    for some fixed $\theta>(1-\gamma-\beta)\lor\frac{1}{2}$ which is independent of $t$.

    \item\label{derivative_second_argument} for every fixed $t \in [0,T]$ the function $s \mapsto f(t,s)$ is absolutely continuous on $[\tau,t]$ with density $\ds f(t,\cdot)$, i.e.
    \[
     f(t,s) = f(t,\tau) +\int_\tau^s \du f(t,u) \ \dd u, \qquad \tau \le s \le t,
    \]
    where $f(t,-\infty) := \lim_{x \to -\infty} f(t,x) = 0$, such that
    \begin{enumerate}
    \item the function $t \mapsto \ds f(t,s)$ is continuous on $(s,\infty)$ for Lebesgue-a.e. $s \in [\tau,\infty)$,
   \item there exist $\eta >0$ and $q >\nicefrac{1}{2}+\nicefrac{5\eta}{2}$  (independent of $t$) such that
    
    \[
     \sup_{t \in [0,T]} \int_{-\infty}^t \left| \ds f(t,s) \right|^{1+\eta} \left(|s|^{q}  \lor 1 \right) \ \dd s < \infty.
     \]
    \end{enumerate}
    \end{enumerate}
\end{definition}

The next lemma provides some useful path and moment properties of the L\'evy-driven Volterra process $M$.

\begin{lemma}\label{l.moments_M}
 Let $f \in \mathcal K$. Then there exists a modification of $M$ (which we still denote by $M$ and which is fixed from now on) such that:
 \begin{enumerate}
  \item $M$ has c\`adl\`ag paths.
  \item The jumps of $M$ fulfil
  \[
   \Delta M(t) = f(t,t) \Delta L(t), \quad t>0.
  \]
  \item Whenever $L(1) \in \mathscr L^p(\mathbb P)$ for some $p \ge 2$ we have
  \[
   \sup_{t \in [0,T]} |M(t)|^p \in \mathscr L^p(\mathbb P).
  \]

 \end{enumerate}

\end{lemma}

\begin{pf}
 In the case $\tau > -\infty$ the assertion follows from Remark~5 in~\cite{BKO_ppm}. In the case $\tau = -\infty$
 it follows from Theorem~8 in~\cite{BKO_ppm} with the choice $\varphi_{q'}(t) = |t|^{q'} \lor 1$ for
 \[
 \frac{1}{2}<q'<\theta \wedge\frac{q-2\eta}{1+\eta}.
 \]
\end{pf}

With regard to the advantages of the present paper over the already existing literature on It\=o formulas for stochastic integrals let us emphasise that in particular fractional L\'evy processes (via the Mandelbrot-Van Ness representation) are included here. Indeed, the following lemma shows that the class $\mathcal K$ contains the kernels

\begin{equation}\label{e.fractional_kernel_definition}
 f_d(t,s)=\frac{1}{\Gamma(d+1)}\left((t-s)^d_+-(-s)^d_+\right)
\end{equation}

for $s,t\in\R$ and a fractional integration parameter $d\in(0,\nicefrac{1}{2})$, where $\Gamma$ denotes the Gamma function. The parameter $d$ is related to the well-known Hurst parameter via $d=H-\nicefrac{1}{2}$.

\begin{lemma}\label{fractional_kernel_in_k}
 The function $f_d:\R^2 \to \R$, defined in \eqref{e.fractional_kernel_definition}, satisfies the assumptions in Definition~\ref{d.1} with~$\tau= -\infty$.
\end{lemma}

The proof of this lemma can be found in the appendix. 

Throughout this paper we use the following definition:

\[
\mathcal A(\R):=\{\xi: \R \to \R: \xi \text{ and } \mathcal F \xi \text{ are in } \mathscr L^1(\dd u)\},
\]

where $\mathcal F\xi$ denotes the Fourier transform of $\xi$. Note that the functions in $\mathcal{A}(\R)$ are continuous and bounded. Furthermore, we use the abbreviation DCT for \emph{dominated convergence theorem}.

\section{$S$-transform and Hitsuda-Skorokhod integrals}\label{s.stochastic_tools}

In this section we make precise the definition of the Hitsuda-Skorokhod integrals which appear in our generalised It\=o formula. 
This definition builds on the injectivity of the Segal-Bargman transform (in short, $S$-transform), which is a tool from white noise analysis.

\subsection{The Segal-Bargmann transform}\label{s.s-transform}

We first introduce a set $\Xi$ by 
\[
\begin{aligned}
 \Xi :=  \text{span}\{ & g: \R_0 \times \R \to \R :  \ g(x,t) = g_1(x) g_2(t) \text{ for two measurable functions such that there} \\
 &\text{exists an }n \in \N \text{ with } \text{supp}(g_1) \subset [-n,-\nicefrac{1}{n}] \cup [\nicefrac{1}{n},n], \ |g_1|\le n\text{ and }  g_2 \in \S \}. 
\end{aligned} 
\]

Here $\S$ is the Schwartz space of rapidly decreasing smooth functions and for any function $h: \R \times \R \to \R$ we define $h^*$  by
\[
 h^*(x,s) := x h(x,s)
\]
for all $x,s \in \R$.

\begin{rem}\label{r.estimate_g}
Let $g \in \Xi$ be given by $ g(x,t) = \sum_{j=1}^N \mu_j g_{1,j}(x)g_{2,j}(t)$. Using the abbreviations 

\[
g_1(x):=  \sum_{j=1}^N \left|\mu_j g_{1,j}(x)\right| \quad \text{and} \quad g_2(t):=  \sum_{j=1}^N \left|g_{2,j}(t)\right|
\]

we see easily that there exists an $n' \in \N$ such that $\text{supp}(g_1) \subset [-n',-\nicefrac{1}{n'}] \cup [\nicefrac{1}{n'},n']$, $|g_1| \le n'$ and that $\sup_{t \in \R}|g_2(t)p(t)|$ is finite for every polynomial $p$ as well as

\[
 |g(x,t)| \le g_1(x) g_2(t)
\]

for every $x \in \R_0$ and all $t\in \R$. We will make use of this simple estimate in our subsequent calculations.
\end{rem}

For every $n\in\N$ let $I_n$ be the $n$-th order multiple L\'evy-It\=o integral
(with respect to $\tilde N$), see e.g. page~665 in \cite{shih08}. For any $g \in
\mathscr L^2( x^2 \nu( \dd x)\times \dd t)$ let $g^{\otimes n}$,
$n\in\N \cup \{0\}$, be the $n$-fold tensor product of $g$ and define a measure $\qg$ on
$(\Omega,\mathscr F)$ by the change of measure

\begin{equation}\label{e.change_of_measure}
\dd \qg =\we(I_1(g)) \dd \mathbb P, 
\end{equation}

where the Radon-Nikod\'ym derivative is the Wick exponential of the random variable $I_1(g)$:
\begin{equation}\label{e.RNderivative}
\we(I_1(g)):=\sum_{n=0}^\infty \frac{I_n(g^{\otimes n})}{n!}.
\end{equation}

 Let us
point out that it follows from
\[
\mathbb E\left(\sum_{n=0}^\infty \frac{I_n(g^{\otimes n})}{n!}\right)=1
\]
that $\qg$ is a signed probability measure. In the following, $\E^{\qg}$
denotes the expectation under $\qg$.

We also mention that  according to \cite{shih08}, Theorem 4.8, we have for $g \in \mathscr L^2\left( x^2 \nu( \dd x)\times \dd t\right)$ with $g^* \in \mathscr L^1\left( \nu( \dd x)\times \dd t\right)$ that
\[
\begin{aligned}
 \we(I_1(g)) &= \exp \left\{ \sigma \int_\R
g(0,t) \   W(\dd t)  - \frac{\sigma^2}{2}\int_{\R_0}g(0,t)^2 \  \dd t    -\int_\R
\int_{\R_0}g^*(x,t)  \ \nu(\dd x)  \ \dd t \right\} \\
 & \qquad \cdot \prod_{t: \Delta L(t) \neq 0} \left( 1+g^*(\Delta L(t),t)
\right),
\end{aligned}
\]

which equals the Dol\'eans-Dade  exponential of $I_1(g)$ at infinity.

Using Proposition 1.4 and formula (10.3) in \cite{nop_malliavin} and the fact that
the Brownian part and the jump part are independent, we infer by applying the
Cauchy-Schwarz inequality that there exists a constant $e_g>0$ (only depending
on $g$) such that
\begin{equation}\label{estimate_changed_measure}
 \E^{\qg}\left(|X|\right) \leq \E\left(|X|^2\right)^{\nicefrac{1}{2}} \cdot
\E\left(\left|\sum_{n=0}^\infty \frac{I_n(g^{\otimes
n})}{n!}\right|^2\right)^{\nicefrac{1}{2}} \leq e_g \cdot
\E\left(|X|^2\right)^{\nicefrac{1}{2}}
\end{equation}
holds for every $X\in\mathscr L^2\left(\mathbb P\right)$.

We are now in the position to define the Segal-Bargmann transform on $\mathscr L^2(\mathbb P)$.

\begin{definition}
 For every $\varphi \in \mathscr L^2(\mathbb P)$ its \emph{Segal-Bargmann transform} (subsequently referred to as \emph{$S$-transform}) $S \varphi$ is given as an integral transform on the set $\mathscr L^2\left( x^2 \nu( \dd x)\times \dd t\right)$ by
\[
 S\varphi(g) := \E^{\qg}(\varphi).
\]

\end{definition}

The following injectivity result for the $S$-transform provides us with a key property for both the definition of Hitsuda-Skorokhod integrals and the proof of the generalised It\=o formula.

\begin{proposition}\label{injectivity_s-transform}
 Let $\varphi, \psi$ be in $\mathscr L^2(\mathbb P)$. If  $S \varphi (g) = S
\psi(g)$ for all $g \in \Xi$, then we have $\varphi = \psi$ $\mathbb P$-almost
surely.
\end{proposition}

\begin{pf}
Noting that
$\Xi$ is dense in $\mathscr L^2\left( x^2 \nu( \dd x)\times \dd t\right)$ 
this result is a direct consequence of (4.1) 
in \cite{shih08}.
\end{pf}

\subsection{Hitsuda-Skorokhod integrals}\label{s.HSi}

Due to the above injectivity property of the $S$-transform we can use it to define 
integration with respect to L\'evy-driven Volterra processes, thereby generalising the approach in \cite{B_sta} and \cite{BM_sccl}.

The motivation for our approach of defining Hitsuda-Skorokhod integrals lies in the fact that under suitable integrability and predictability assumptions on the integrand they reduce to the well known stochastic integrals with respect to semimartingales and random measures, respectively.

 \begin{definition}\label{M_diamond}
 
Let $\mathcal B \subset[0, \infty)$ be a Borel set. Suppose the mapping $t \mapsto S(M(t))(g)$ is differentiable for every $g \in \Xi$, $t \in \mathcal B$ , and let $X:\mathcal B\times \Omega \rightarrow \R$ be a stochastic process such that $X(t)$ is square-integrable for a.e. $t \in \mathcal B$. The process $X$ is said to have a \emph{Hitsuda-Skorokhod integral} with respect
to $M$ if there is a $\Phi \in \mathscr L^2(\mathbb P)$ such
that 
  \[
    S\Phi(g) = \int_\mathcal B S(X(t))(g) \ddt S(M(t))(g) \ \dd t
  \]
holds for all $g \in \Xi$. As the $S$-transform is injective, $\Phi$ is unique and we write
  \[
   \Phi = \int_\mathcal B X(t) \ M^\diamond(\dd t).
  \]
\end{definition}

\begin{rem}
 A different approach of defining Skorokhod integrals with respect to fractional L\'evy processes via white noise analysis can be found in~\cite{LD_sif}. If the fractional L\'evy process is of finite $p$-variation, stochastic integrals with respect to this process can be defined pathwise as an improper Riemann-Stieltjes integral and have been considered in \cite{FK_flou}. In the special case that $M$ is  a L\'evy process, it can be shown with the techniques in \cite{nop_malliavin} that our definition of Hitsuda-Skorokhod integrals coincides with the definition of Skorokhod integrals via the chaos decomposition. A related approach to Skorokhod integrals with respect to Poisson-driven Volterra processes via Malliavin calculus is provided in \cite{DS_ac}.
\end{rem}

The following technical lemma will prove useful later on. The proof of this lemma is provided in the appendix.

\begin{lemma}\label{l.continuity_of_integral}

Let $F:\R^2 \to \C$ with $\text{supp }\, F \subset [\tau, \infty)^2$ for some $\tau \in [-\infty, 0]$, $F(t,\cdot)\in\mathscr L^1(\dd s)$ for every $t\in[0,T]$ and let  $\beta,\gamma\in(0,1)$ with $\beta+\gamma<1$. We define

\[
 \text{I}_F(t) := \int_{-\infty}^t F(t,s) \ \dd s.
\]

 \begin{enumerate}[{\rm(i)}]

\item\label{F_i_conditions_2}  Let the following set of conditions be satisfied:
 \begin{enumerate}[{\rm a)}]
   \item\label{F_continuity_kernel_ii} For Lebesgue-a.e. $s \in \R$ the map $t \mapsto F(t,s)$ is continuous on the set $[0,T]\setminus \{s\}$.
  \item\label{F_integrability_derivative_ii} For every $t \in [0,T]$ there exist an $\varepsilon >0$ and a constant $\tilde C>0$ such that
    \begin{equation}\label{e.F_integral_bounded}
     \int_{-\infty}^{t-2\varepsilon} \sup_{r \in [0 \lor (t-\varepsilon),\, (t+\varepsilon) \land T]} \left|  F(r,s) \right| \ \dd s < \infty
    \end{equation}   
and  
    \begin{equation}\label{e.F_estimate}
     \left| F(r,s)\right| \le \tilde C|s|^{-\beta} |r-s|^{-\gamma}
    \end{equation}
   for all $r \in [(t-\varepsilon) \lor 0, \, (t+\varepsilon) \land T)]$ and $s \in [t-2 \varepsilon, r)$.
 \end{enumerate}
 Then the function $I_F$ is continuous on $[0,T]$.
 
  \item\label{F_i_conditions} Let the following set of conditions be satisfied:
 
 \begin{enumerate}[{\rm a)}]
  \item\label{F_continuity_kernel} The mapping $(s,t) \mapsto F(t,s)$ is continuous on the set $\{(t,s) \in \R^2: \  \tau \le s \le t \}$.
  \item\label{F_differentiable} For Lebesgue-a.e. $s \in \R$ the map $t \mapsto F(t,s)$ is continuously differentiable on $[0,T]\setminus \{s\}$.
  \item\label{F_integrability_derivative} For every $t \in [0,T]$ there exist an $\varepsilon >0$ and a constant $\tilde C>0$ such that
    \begin{equation}\label{e.diff_F_integral_bounded}
     \int_{-\infty}^{t-2\varepsilon} \sup_{r \in [(t-\varepsilon) \lor 0, \, (t+\varepsilon) \land T]} \left| \dr F(r,s) \right| \ \dd s < \infty
    \end{equation}
  
    and
  
    \begin{equation}\label{e.diff_F_estimate}
     \left|\dr F(r,s)\right| \le \tilde C|s|^{-\beta} |r-s|^{-\gamma}
    \end{equation}

     for all $r \in [(t-\varepsilon) \lor 0, \, (t+\varepsilon) \land T)]$ and $s \in [t-2 \varepsilon, r)$.
     
 \end{enumerate}
Then the function $I_F$ is continuously differentiable on $[0,T]$ with derivative

\[
 I_F'(t) = F(t,t) + \int_{-\infty}^t  \dt F(t,s) \ \dd s, \quad t \in [0,T].
\]

\end{enumerate}

\end{lemma}

We next derive the explicit form of the derivative of the $S$-transform of $M(t)$. This result particularly yields that the differentiability condition on the mapping $t \mapsto  S(M(t))(g)$ in Definition~\ref{M_diamond} is fulfilled for kernel functions $f \in \mathcal K$.

\begin{lemma}\label{derivative_S(M(t))}
 For all $f \in \mathcal K$ and $g\in\Xi$ the mapping $t \mapsto S(M(t))(g) $ is continuously differentiable on the set~$[0,T]$ with derivative 
\[
 \begin{aligned}
  \frac{\dd }{\dd t} S(M(t))(g) = &\sigma \left(f(t,t) g(0,t)+
\int_{-\infty}^t  \dt f(t,s) g(0,s) \ \dd s\right) \\ &+f(t,t) \int_{\R_0}x
g^*(x,t) \ \nu (\dd x)   +  \int_{-\infty}^t \int_{\R_0} \dt f(t,s) x g^*(x,s)
\ \nu(\dd x) \ \dd s.
 \end{aligned}
\]
\end{lemma}
Note that the finiteness of the integrals appearing in the above lemma is
guaranteed by Definition~\ref{d.1}\eqref{behaviour_kernel}.

\begin{pf}
By the isometry of L\'evy-It\=o integrals  we
obtain
\begin{equation}\label{SMt}
  S(M(t))(g) =  \sigma  \int_{-\infty}^t  f(t,s) g(0,s)  \ \dd s +
\int_{-\infty}^t \int_{\R_0}  f(t,s) x g^*(x,s) \  \nu(\dd x) \ \dd s =: {\rm I}(t)+ {\rm II}(t),
\end{equation}
cf. e.g. Section~3.1 of \cite{B_sta} and  Example~3.6 in \cite{BM_sccl}.

We now want to apply Lemma~\ref{l.continuity_of_integral} to ${\rm I}(t)$ with $F_{\rm I}(t,s):= f(t,s)g(0,s)$. It is easy to check that conditions~a) and~b) of Lemma~\ref{l.continuity_of_integral}\eqref{F_i_conditions} are satisfied. Since $\sup_{s \in \R} |g(0,s)| < \infty$ (cf. Remark~\ref{r.estimate_g}), we deduce from \eqref{e.decay_of_dtf} that~\eqref{e.diff_F_estimate} is  fulfilled for every $r>s>-\infty$. Moreover, by using \eqref{e.decay_of_dtf} and the rapid decrease of $g(0,s)$ we get for $t \in [0,T]$ and arbitrary $\epsilon >0$

\[
\begin{aligned}
\int_{-\infty}^{t-2 \epsilon} \sup_{r \in [(t-\epsilon) \lor 0, \, (t+\epsilon) \land T]} \left|\dr F_{\rm I}(r,s)\right|  \ \dd s &=  \int_{-\infty}^{t-2 \epsilon}  \sup_{r \in [(t-\epsilon) \lor 0, \, (t+\epsilon) \land T]} \left| \dr f(r,s)   g(0,s) \right| \ \dd s \\ 
&\leq  C_0 \int_{-\infty}^{t-2 \epsilon}   |s|^{-\beta} \sup_{r \in [(t-\epsilon) \lor 0, \, (t+\epsilon) \land T]}  |r-s|^{-\gamma} |g(0,s)| \  \dd s \\
&\leq  C_0 \int_{-\infty}^{t-2 \epsilon}   |s|^{-\beta}|t-\epsilon-s|^{-\gamma}  |g(0,s)| \  \dd s \\
  &< \infty, 
\end{aligned}
\]

where $C_0$, $\beta$ and $\gamma$ are given by \eqref{e.decay_of_dtf}. This shows that~\eqref{e.diff_F_integral_bounded} holds and thus Lemma~\ref{l.continuity_of_integral} is applicable, which results in $t \mapsto {\rm I}(t)$ being continuously differentiable on $[0,T]$.

To deal with II$(t)$ we recall from Remark~\ref{r.estimate_g} that every $g \in \Xi$ can be written in the form $ g(x,t) = \sum_{j=1}^N \mu_j g_{1,j}(x)g_{2,j}(t)$. We set

\[
 F_{\rm II}(t,s) := \int_{\R_0} f(t,s) x g^*(x,s) \ \nu(\dd x) = \sum_{j=1}^N \mu_j \left( \int_{\R_0}  x^2 g_{1,j}(x) \ \nu(\dd x) \right) f(t,s) g_{2,j}(s).
\]

where $g_{2,j} \in \S$ for all $1 \le j \le N$. Hence, each of the summands in $F_{\rm II}$ is of the same form as $F_{\rm I}$. Consequently, $t \mapsto {\rm II}(t)$ is continuously differentiable. Therefore, in view of \eqref{SMt} the mapping $t \mapsto  S(M(t))(g)$ is continuously differentiable on $[0,T]$.
\end{pf}

We proceed by introducing a Hitsuda-Skorokhod integral with respect to an appropriate random measure that will enable us to establish a connection between the Hitsuda-Skorokhod integral with respect to the L\'evy-driven Volterra process $M$ and the classical integral with respect to the underlying L\'evy process $L$ (see Theorem~\ref{r.connection_m_diamond_lamnda_diamond} below).

\begin{definition}\label{Lambda-diamond}
Let $\mathcal B \subset \R$ be a Borel set and $X:\R \times \mathcal B \times \Omega \rightarrow \R$ be a random field such that $X(x,t)  \in \L^2(\mathbb P)$ for $\nu(\dd x) \otimes \dd t$-a.e.~$(x,t)$. The \emph{Hitsuda-Skorokhod integral} of $X$ with respect to the random measure 
\[
\Lambda (\dd x, \dd t) = x\tilde N(\dd x, \dd t)+ \sigma \delta_0( \dd x) \otimes W (\dd t),
\]
where $\delta_0$ denotes the Dirac measure in $0$, is said to exist in $\L^2(\mathbb P)$, if there is a random variable $\Phi \in
\L^2(\mathbb P)$ that satisfies

   \[
   \begin{aligned}
          S\Phi(g) &= \int_\mathcal B \int_{\R_0} S(X(x,t))(g) g^*(x,t) \ x \ \nu(\dd x) \ \dd t +  \sigma \int_\mathcal B S(X(0,t))(g) \ g(0,t) \ \dd t
   \end{aligned}
   \]
  for all $g \in \Xi$.  In this case, by Proposition~\ref{injectivity_s-transform} the random variable $\Phi$ is unique and we write
\[
 \Phi = \int_\mathcal B \int_{\R} X(x,t) \ \Lambda^\diamond (\dd x,\dd t). 
\]
\end{definition}

\begin{rem}\label{special_cases}
 Let $X:\R \times [0,T] \times \Omega \rightarrow \R$ be a predictable random field as in Definition~\ref{Lambda-diamond} with $\mathcal B= [0,T]$.
 \begin{enumerate}
  \item\label{B-diamond} Assume that $\sigma >0$ and let $X$ be given by
  \[
   X(x, \cdot) = \begin{cases}
                  \frac{1}{\sigma} Y(\cdot), &x =0 \\
                  0, &x \ne 0
                 \end{cases}
  \]
  for some stochastic process $Y: {[0,T]} \times \Omega \to \R$. Since $X$ is predictable, we infer from Theorem~3.1 in \cite{B_sta} that the Hitsuda-Skorokhod integral $\int_0^T \int_{\R} X(x,t) \ \Lambda^\diamond(dx,dt)$ exists and satisfies
  \[
     \int_0^T \int_{\R} X(x,t) \ \Lambda^\diamond (\dd x,\dd t) =   \int_0^T  Y(t) \ W (\dd t),
  \]
    where the last integral is the classical stochastic integral with respect to the Brownian motion $W$.
    Note that it follows from the calculations in Section~3.1 of \cite{B_sta} that
    \[
	 g(0,t) = \ddt \int_0^t g(0,s) \ \dd s  =\ddt S(W(t))(g)
    \]
   and hence we have
   \[
    S\left( \int_0^T \int_{\R} X(x,t) \ \Lambda^\diamond (\dd x,\dd t) \right) (g)= \int_0^T S(Y(t))(g) \ \ddt S(W(t))(g) \ \dd t.   
   \]

   \item If $X$ fulfils $X(0, \cdot) \equiv 0$, then it follows from Theorem~3.5 in \cite{BM_sccl} that
  \[
     \int_0^T \int_{\R} X(x,t) \ \Lambda^\diamond (\dd x,\dd t) =  \int_0^T \int_{\R_0} x X(x,t) \ \tilde{N} (\dd x, \dd t),
  \]
  where the last integral is the classical stochastic integral with respect to the compensated Poisson jump measure $\tilde N$.
\item According to 1. and 2. we have
\[
      \int_0^T \int_{\R} X(x,t) \ \Lambda^\diamond (\dd x,\dd t) =  \int_0^T \int_{\R_0} x X(x,t)  \ \tilde N (\dd x,\dd t)+ \sigma \int_0^T X(0,t) \ W (\dd t).
     \]
 \end{enumerate}
\end{rem}

Now we are in the position to state the connection between the Hitsuda-Skorokhod integral with respect to $M$ and the 
It\=o integral with respect to the driving L\'evy process. 

\begin{theorem}\label{r.connection_m_diamond_lamnda_diamond}
Suppose that $X$ is a predictable process such that

\begin{equation}\label{e.second_running_maximum}
\E\left( \sup_{t \in [0,T]}|X(t)|^2\right) < \infty
\end{equation}

and let $f\in\mathcal K$. Then

 \begin{equation}\label{e.existence_M-diamond}
\int_0^T X(t) \ M^\diamond (\dd t) 
\end{equation}
exists, if and only if 
 \begin{equation}\label{e.existence_Lambda-diamond}
\int_{-\infty}^T \int_{\R} \int_{0 \lor s}^{T} \dt f(t,s) X(t) \ \dd t \ \Lambda^\diamond (\dd x,\dd s) 
\end{equation}
exists. In this case 
 \begin{equation}\label{e.connection_m_diamond_lambnda_diamond}
    \int_0^T X(t) \ M^\diamond (\dd t) = \int_0^T  f(t,t) X(t) \ L(\dd t) + \int_{-\infty}^T \int_{\R} \int_{0 \lor s}^{T} \dt f(t,s) X(t) \ \dd t \ \Lambda^\diamond (\dd x,\dd s) .
    \end{equation}
 
\end{theorem}

 {\bf Proof}

 Note that $\int_0^T  f(t,t) X(t) \ L(\dd t) $ exists in $\L^2(\mathbb P)$ because of the continuity of $t \mapsto f(t,t)$ and \eqref{e.second_running_maximum}. By the previous remark 
its $S$-transform is given by
 
 \[
 \begin{aligned}
   S&\left( \int_0^T  f(t,t) X(t) \ L(\dd t)\right)(g) \\
   &= \sigma \int_0^T  S(X(t))(g) f(t,t) g(0,t)  \ \dd t + \int_0^T \int_{\R_0} S(X(t))(g) f(t,t)  x g^*(x,t) \ \nu(\dd x) \ \dd t.
 \end{aligned} 
 \]
 
 Using assumption \eqref{e.second_running_maximum} and \eqref{e.decay_of_dtf} we obtain the estimate

 \[
  \E \left( \left( \int_{0 \lor s}^{T} \dt f(t,s) X(t) \ \dd t \right)^2\right) \le \E\left( \sup_{t \in [0,T]}|X(t)|^2\right)   \left( \int_{0 \lor s}^{T} \dt f(t,s)  \ \dd t \right)^2 < \infty.
 \]

 Thus, we can apply Definition~\ref{Lambda-diamond} to $\int_{0 \lor s}^{T} \dt f(t,s) X(t) \ \dd t$. Assuming the existence of~\eqref{e.existence_Lambda-diamond}, we therefore infer from Definition~\ref{M_diamond}, Lemma~\ref{derivative_S(M(t))} and Fubini's theorem that~\eqref{e.existence_M-diamond} exists and satisfies~\eqref{e.connection_m_diamond_lambnda_diamond}. Analogous arguments yield the converse implication.
 \hfill{$\square$}

Considering the jump measure $N$ instead of the compensated jump measure $\tilde N$ naturally leads to the following definition by adding the $S$-transform of the integral with respect to the compensator. 

\begin{definition}\label{N-diamond}
Let $\mathcal B \subset \R$ be a Borel set and $X:\R_0 \times \mathcal B \times \Omega \rightarrow \R$ be a random field such that $X(x,t)  \in \L^2(\mathbb P)$ for $\nu(\dd x) \otimes \dd t$-a.e.~$(x,t)$. The \emph{Hitsuda-Skorokhod integral} of $X$ with respect to the jump measure
$ N(\dd x, \dd t)$ is said to exist in $\L^2(\mathbb P)$, if there is a (unique) random variable $\Phi \in
\L^2(\mathbb P)$ that satisfies
   \[
      S\Phi(g) = \int_\mathcal B \int_{\R_0} S(X(x,t))(g) (1+g^*(x,t)) \ \nu(\dd x) \
\dd t
   \]
  for all $g \in \Xi$. We write 
  \[
  \Phi = \int_\mathcal B \int_{\R_0} X(x,t) \ N^\diamond (dx,dt).
  \]
\end{definition}

\section{Generalised It\={o} formulas}\label{s.Ito-formula}

This section is devoted to formulating precisely and to proving our generalised It\=o formula. 

The following theorem is the main result of this paper:

\begin{theorem}\label{main_result_L2}
Let $f \in \mathcal K$ and $G\in C^2(\R)$.
Additionally, assume that one of the following assumptions is fulfilled:
\begin{enumerate}[(i)]
 \item $\sigma>0$, $G, G'$ and $G''$ are of polynomial growth with degree ${q} \ge 0$, that is

\[
 |G^{(l)}(x)| \le C_{pol} (1+|x|^{q}) \quad \text{ for every $x \in \R$ and $l=0,1,2$}
\]

with a constant $C_{pol}>0$, and 

\[
L(1) \in \mathscr L^{p_{q}}(\mathbb P)
\]

for $p_{q} = 4 \lor (2{q}+2)$;
 \item $G,G',G''\in\mathcal A(\R)$ and $L(1) \in \mathscr L^{4}(\mathbb P)$.
  \end{enumerate}

Then the following generalised It\=o formula

\begin{equation}\label{ito_formula}
  \begin{aligned}
G(M(T)) &= G(0) + \frac{\sigma^2}{2} \int_{0}^T G''(M(t)) \left( \ddt
\int_{-\infty}^t f(t,s)^2  \dd s \right) \dd t\\
&\quad+ \sum_{0<t\leq T} \left[ G(M(t))- G(M(t-))- G'(M(t-)) \Delta M(t)
\right]
\\
&\quad+ \int_{-\infty}^T \int_{\R} \int_{0 \lor s}^T G'\left(M(t) + xf(t,s)\right)  \dt f(t,s)  \ \dd t \ \Lambda^\diamond(\dd x, \dd s)  \\
    &\quad+  \int_{-\infty}^T \int_{\R_0} \int_{0 \lor s}^T \left( G'\left(M(t) +
    xf(t,s)\right) -G'(M(t)) \right) x \dt f(t,s) \ \dd t \ \nu(\dd x)\  \dd s  \\
&\quad+ \int_{0}^T G'(M(t-)) f(t,t) \ L (\dd t)
 \end{aligned}
\end{equation}

 holds $\mathbb P$-almost surely. In particular, all the terms in~\eqref{ito_formula} exist in $\L^2(\mathbb P)$. Moreover, the generalised It\=o formula is valid in the form of (\ref{ito_formula_L2}), provided all terms there exist in $\L^2(\mathbb P)$.
\end{theorem}

\begin{rem}
We point out that formula (\ref{ito_formula}) also holds in case (ii) under the weaker assumption that 
$L(1) \in \mathscr L^{2}(\mathbb P)$, provided that $M$ is continuous, i.e. $f(t,t)\equiv 0$.
\end{rem}

In order to prove Theorem~\ref{main_result_L2} we start with a heuristic argumentation that gives a rough outline of the steps of our proof and motivates the auxiliary results that we shall prove below. Suppose the fundamental theorem of calculus enables us to write

\begin{equation}\label{e.motivation_1}
  S(G(M(T)))(g) = G(0)+ \int_0^T \ddt  S(G(M(t)))(g)\  \dd t.
\end{equation}

Subsequently, by making use of the Fourier inversion theorem in the spirit of \cite{LS_ito} we would obtain
\begin{equation}\label{e.motivation_2}
    S(G(M(t)))(g) = \E^{\qg} \bigl(G(M(t)) \bigr) = \frac{1}{\sqrt{2 \pi}} \int_\R  (\F G)(u)  \E^{\qg} \bigl(e^{iuM(t)} \bigr) \ \dd u.
\end{equation}

Differentiating the right-hand side of~\eqref{e.motivation_2} with respect to $t$, using some standard manipulations of the Fourier transform, plugging the resulting formula for $\ddt  S(G(M(t)))(g)$ into \eqref{e.motivation_1} and using the injectivity of the $S$-transform would then give an explicit expression for $G(M(T))$ leading to a generalised It\=o formula. 

Our approach to prove Theorem~\ref{main_result_L2} is based on several auxiliary results. More precisely, following the above motivation we derive explicit expressions for the characteristic function $\E^{\qg}\left(e^{iuM(t)}\right)$ under $\qg$ as well as its derivative $\dt \E^{\qg} \left(e^{iuM(t)} \right)$ in Proposition~\ref{eiumt} and Lemma~\ref{derivative_S(exp(iuM(t)))}. 

In (the proof of) Proposition~\ref{e.derivative_SGM} we will show that the integral appearing in~\eqref{e.motivation_2} is well-defined and that the mapping $t \mapsto S(G(M(t)))(g)$ is indeed differentiable. We then complete the proof of Theorem~\ref{main_result_L2} by using the explicit expression for $\dt \E^{\qg} \left(e^{iuM(t)} \right)$ and the injectivity of the $S$-transform.

Let us now follow our approach by providing the characteristic function of $M(t)$. The following result was obtained in Proposition~2.7 of \cite{RR89}: 
\begin{lemma}\label{eiumt_0}
 For every $f \in \mathcal K$ and $t\ge0$ we have
\[
 \begin{aligned}
  &\E\left( e^{iuM(t)} \right)  \\
& \qquad = \exp \left(  -\frac{\sigma^2 u^2}{2} \int_{-\infty}^t f(t,s)^2 \  \dd s  +
\int_{-\infty}^t \int_{\R_0} \left( e^{iuxf(t,s)} -1 -iux f(t,s) \right) \ \nu (\dd x) \  \dd s \right).
 \end{aligned}
\]
\end{lemma}

The following proposition is concerned with  the characteristic function of $M(t)$ under the signed
measure $\qg$.

\begin{proposition}\label{eiumt} 
 Let $f \in \mathcal K$ and $g\in\Xi$. Then
\begin{equation}\label{S(exp(iuM(t)))(g)}
 \begin{aligned}
  &\E^{\qg}\left( e^{iuM(t)} \right)  \\
&= \exp \left( iu  \int_{-\infty}^t \sigma^2 f(t,s) g(0,s)  \ \dd s
+iu\int_{-\infty}^t  \int_{\R_0} x f(t,s) g^*(x,s) \  \nu(\dd x) \  \dd s  \right. \\
& \quad  \left. -\frac{\sigma^2 u^2}{2} \int_{-\infty}^t f(t,s)^2 \  \dd s  +
\int_{-\infty}^t \int_{\R_0} \left( e^{iuxf(t,s)} -1 -iux f(t,s) \right)(1+
g^*(x,s))  \ \nu (\dd x) \  \dd s \right).
 \end{aligned}
\end{equation}  
\end{proposition}

\begin{pf}
Note that
\begin{equation}\label{e.eiumt.1}
M(t)=I_1(f(t,\cdot)).
\end{equation}
Approximating $f(t,\cdot)$ by the sequence of functions $(f_n(t,x,s))_{n \in \N}$  defined via $$f_n(t,x,s):= \mathds 1_{[-n,n]}(x) \mathds 1_{[-n,t]}(s) f(t,s)$$
and using Theorem~5.3 in \cite{LS_sbt} we see that $\mathbb P$-a.s. the equation
\[
 e^{I_1(iuf_n(t,\cdot))} = \E\left( e^{I_1(iuf_n(t,\cdot))} \right) \we \left( I_1\left(k_{t,n} \right) \right) 
\]
holds with
\[
\begin{aligned}
  k_{t,n}(x,s) &:=  \mathds1_{\{0\}}(x) iuf_n(t,x,s)+ \mathds1_{\R_0}(x) \frac{e^{iuxf_n(t,x,s)}-1}{x} \\
    &=\mathds 1_{[-n,t]}(s) \mathds 1_{[-n,n]}(x) k_t(x,s),
\end{aligned}
\]
where $k_t$ is given by $k_t(x,s):= \mathds1_{\{0\}}(x) iuf(t,s)+ \mathds1_{\R_0}(x) \frac{e^{iuxf(t,s)}-1}{x}$. 
Note that the results in \cite{LS_sbt} are applicable, because by construction of $f_n$ we can switch from the original L\'evy process with L\'evy measure $\nu$ to an auxiliary one 
with L\'evy measure $\nu_n(A)=\nu(A\cap [-n,n])$, for which the moment conditions assumed in \cite{LS_sbt} are satisfied.

In view of (\ref{e.change_of_measure}) this yields
\[
 \E^{\qg}\left( e^{I_1(iuf_n(t,\cdot))} \right) =\E\left( e^{I_1(iuf_n(t,\cdot))} \right)\E\Bigl(\we(I_1(g)) \we(I_1(k_{t,n})) \Bigr).
\]
By the isometry for multiple L\'evy-It\=o integrals and (\ref{e.RNderivative}) we obtain
\begin{equation}\label{e.eiumt.2}
\begin{aligned}
  &\E^{\qg}  \left( e^{I_1(iuf_n(t,\cdot))} \right) \\
  &=\E\left( e^{I_1(iuf_n(t,\cdot))} \right) \exp\left( \int_{-n}^t \int_{-n}^n g(x,s) k_{t}(x,s) \ \tilde{\nu}(\dd x, \dd s) \right),
\end{aligned}
\end{equation}
where the measure $\tilde{\nu}$ is defined as $\tilde{\nu}(\dd x, \dd s)= x^2 \nu(\dd x) \times \dd s$ if $x \ne0$ and $\tilde{\nu}(\dd x, \dd s)= \sigma^2 \dd s$ if $x =0$. Resorting to (\ref{e.eiumt.1}) and (\ref{e.eiumt.2}), taking the limit as $n \to \infty$ and using DCT results in
\[
  \E^{\qg} \left( e^{iuM(t)} \right) =\E\left( e^{iuM(t)} \right) \exp\left( \int_{-\infty}^t \int_\R g(x,s) k_t(x,s) \ \tilde{\nu}(\dd x, \dd s) \right).
\]

Plugging in the formula for $\E(e^{iuM(t)})$ from  Lemma~\ref{eiumt_0}  shows that~\eqref{eiumt} holds.
\end{pf}

In the spirit of Lemma~\ref{derivative_S(M(t))} we now derive a formula for the
derivative of the $S$-transform of $e^{iuM(t)}$.

\begin{lemma}\label{derivative_S(exp(iuM(t)))}
 Let $f \in \mathcal K$ and $g\in\Xi$. Then the map $t \mapsto \E^{\qg}\left( e^{iuM(t)}\right)$ is continuously differentiable on $[0,T]$ with derivative
\[
 \begin{aligned}
& \dt\E^{\qg}\left( e^{iuM(t)} \right) &&
\\
& =  \E^{\qg}\left( e^{iuM(t)} \right) \cdot \Biggl[ && 
iu \frac{\dd }{\dd t}  S(M(t))(g)
\\
&&&-\frac{\sigma^2u^2}{2}\left(f(t,t)^2+ 2\int_{-\infty}^t \dt f(t,s)\cdot f(t,s) \ 
\dd s\right) 
\\
&&&+\int_{\R_0} \left( e^{iuxf(t,t)} -1 -iux f(t,t) \right)(1+ g^*(x,t)) \  \nu
(\dd x) 
\\
&&&+\int_{-\infty}^t \int_{\R_0} \left(iux  \dt f(t,s) \left( e^{iuxf(t,s)} -1
\right)(1+ g^*(x,s)) \right) \  \nu (\dd x)  \ \dd s  \Biggr]. 
 \end{aligned}
\]
\end{lemma}

\begin{pf}
By the differentiability of the exponential function we only have to prove the
differentiability of the terms in the exponential of \eqref{S(exp(iuM(t)))(g)}.
The first two of these summands are easily identified as the terms that occur in \eqref{SMt} and are already treated in Lemma \ref{derivative_S(M(t))}. 

To deal with the fourth summand in the exponential of \eqref{S(exp(iuM(t)))(g)} we define

\[
F(t,s) := \int_{\R_0} \left( e^{iuxf(t,s)} -1 -iux f(t,s) \right)(1+g^*(x,s))  \ \nu (\dd x).
\]

To check for the continuity condition in Lemma~\ref{l.continuity_of_integral}\eqref{F_i_conditions}a) we choose a sequence $(t_n, s_n)_{n \in \N}$ with $(t_n, s_n) \to (t,s)$ as $n \to \infty$. Without loss of generality we may assume $|t-t_n| \le 1$ and $|s-s_n| \le 1$ for all $n \in \N$. In order to apply the DCT to the expression

\begin{equation}\label{F_tn_sn}
  F(t_n,s_n) = \int_{\R_0} \left( e^{iuxf(t_n,s_n)} -1 -iux f(t_n,s_n) \right)(1+g^*(x,s_n))  \ \nu (\dd x)
\end{equation}

we define 

\[
D := \sup_{z \in \R_0} \left| \frac{e^{iz}-1-iz}{z^2}\right| < \infty
\]

and write by using Remark~\ref{r.estimate_g} 

\begin{equation}\label{e.bound_for_continuity}
 \begin{aligned}
   \int_{\R_0} & \sup_{n \in \N} \left|\left( e^{iuxf(t_n,s_n)} -1 -iux f(t_n,s_n) \right)(1+g^*(x,s_n)) \right|  \ \nu (\dd x) \\
   &\le D \int_{\R_0} \sup_{n \in \N} \left(\left(u^2 x^2 f(t_n,s_n)^2 \right)(1+|g^*(x,s_n)|) \right)  \ \nu (\dd x) \\
   &\le D u^2 \int_{\R_0} \sup_{h_1,h_2 \in [-1,1]} \left(\left(x^2 f(t +h_1,s + h_2)^2 \right)(1+\|g^*(x,s+ h_2)|) \right)  \ \nu (\dd x) \\
   &\le Du^2 \sup_{h_1,h_2 \in [-1,1]} \left(f(t +h_1,s + h_2)^2 \right) \int_{\R_0} x^2 \Bigl(1+ g_1(x) \sup_{h_2 \in [-1,1]} g_2(s+ h_2)\Bigr) \ \nu (\dd x) \\
   &< \infty,
 \end{aligned}
\end{equation}

where the finiteness follows from the continuity of $f$ and the fact that the expression in the brackets in the last integral is bounded by Remark~\ref{r.estimate_g}. In view of \eqref{e.bound_for_continuity} the pointwise convergence in $n$ of the integrand of \eqref{F_tn_sn} shows the continuity of the function $F$.

Using again Remark~\ref{r.estimate_g} and (\ref{e.decay_of_dtf}) we get for $t\in [0,T]$, $s<t$ with  $s \ne 0$ and an arbitrary $\epsilon\in(0,t-s)$ the following chain of estimates that will be useful below:

\begin{equation}\label{e.diff_F}
\begin{aligned}
&\int_{\R_0} \sup_{r\in[0 \lor(t-\epsilon),\, (t+\epsilon) \land T]} \left|ux\dr f(r,s) \left( e^{iux f(r,s)} -1\right)\right|  |1+ g^*(x,s)| \ \nu(\dd x) \\
&\leq   2u^2\sup_{y\in\R_0}\left|\frac{e^{iy}-1}{y}\right| \int_{\R_0} x^2 (1 \lor g_1^*(x)) \ \nu(\dd x) \\
&\hspace{4cm}\sup_{r\in[0 \lor(t-\epsilon),\, (t+\epsilon) \land T]}|f(r,s)|\cdot C_0|s|^{-\beta}|t-\epsilon-s|^{-\gamma} (1 \lor g_2(s)) \\
\\
&< \infty, 
\end{aligned}
\end{equation}

where $C_0,\beta$ and $\gamma$  are given by~\eqref{e.decay_of_dtf}.

In order to show the differentiability of the function $F$, we use the mean value theorem for complex valued functions to obtain the estimate

\[
\begin{aligned}
& \int_{\R_0} \sup_{h\in(-\epsilon, \epsilon)\cap (-t, T-t)\setminus \{0\}}\left|\frac{e^{iuxf(t+h,s)}-e^{iuxf(t,s)}-
iux(f(t+h,s)-f(t,s))}{h}(1+ g^*(x,s))\right|  \ \nu(\dd x) \\
&\le \int_{\R_0}\sup_{r\in[0 \lor (t-\epsilon),\, (t+\epsilon) \land T]}\left|ux\dr f(r,s) \left( e^{iux f(r,s)} -1\right)\right|  |1+ g^*(x,s)| \ \nu(\dd x), \\
\end{aligned}
 \]

which is finite by \eqref{e.diff_F} for $t \in [0,T]$, $s<t$ with $s\ne0$ and $\epsilon\in(0,t-s)$.  Therefore an application of the DCT yields the differentiability of $F$ with respect to its first variable and

\begin{equation}\label{e.diff_F_2}
  \dt F(t,s) = iu \int_{\R_0}x \dt f(t,s) \left( e^{iux f(t,s)} -1\right) (1+ g^*(x,s)) \ \nu(\dd x).
\end{equation}

By another application of the DCT in view of \eqref{e.diff_F} we also get the continuity of the derivative $\dt F(t,s)$ for $s \ne 0$ and $t \in [0,T]\setminus \{s\}$. 

To check~\eqref{e.diff_F_integral_bounded} in Lemma~\ref{l.continuity_of_integral}\eqref{F_i_conditions} for $F$ we fix $t \in [0,T]$ and choose $C_0,\beta,\gamma, \theta$ and $\epsilon \in \left(0, \nicefrac{(t+1)}{2} \right)$ such that~\eqref{e.decay_of_dtf} and~\eqref{e.decay_of_f} hold.  Integrating~\eqref{e.diff_F} with respect to $s$ from $-\infty$ to $t-2\epsilon$ thus results in

\[
\begin{aligned}
& \int_{-\infty}^{t-2\epsilon} \int_{\R_0}\sup_{r\in[0 \lor(t-\epsilon),\, (t+\epsilon) \land T]}\left|ux\dr f(r,s) \left( e^{iux f(r,s)} -1\right)\right|  |1+ g^*(x,s)| \ \nu(\dd x) \  \dd s 
\\
&\leq   2u^2\sup_{y\in\R_0}\left|\frac{e^{iy}-1}{y}\right| \int_{\R_0} x^2 (1 \lor |g_1^*(x)|) \ \nu(\dd x) \\
&\qquad \qquad\int_{-\infty}^{t-2\epsilon}\sup_{r\in[0 \lor(t-\epsilon),\, (t+\epsilon) \land T]}\left(|f(r,s)|\cdot C_0|s|^{-\beta} |r-s|^{-\gamma} \right)(1 \lor g_2(s))  \  \dd s \\
\\
&\leq   2u^2 C_0\sup_{y\in\R_0}\left|\frac{e^{iy}-1}{y} \right| \int_{\R_0} x^2 (1 \lor |g_1^*(x)|)\ \nu(\dd x)  \\
&\qquad  \cdot \Biggl( \sup_{v \in (-\infty,-1]} \sup_{r \in [0 \lor(t-\epsilon),\, (t+\epsilon) \land T]} \left( |f(r,v)||v|^{\theta}\right) \int_{-\infty}^{-1} |s|^{-(\theta+\beta)} |t-\epsilon -s|^{-\gamma} (1 \lor g_2(s)) \  \dd s \\
& \qquad \qquad +\sup_{v \in [-1,t-2\epsilon]}\sup_{r\in[0 \lor(t-\epsilon),\, (t+\epsilon) \land T]}|f(r,v)| \int_{-1}^{t-2\epsilon} |s|^{-\beta} |t-\epsilon-s|^{-\gamma} (1 \lor g_2(s)) \  \dd s \bigr) \Biggr)   \\
&< \infty.
\end{aligned}
\]

Note that the finiteness follows from Remark~\ref{r.estimate_g}, \eqref{e.decay_of_dtf}, \eqref{e.decay_of_f},  and Definition~\ref{d.1}\eqref{continuity_of_f}. Hence,~\eqref{e.diff_F_integral_bounded} is fulfilled.

To check~\eqref{e.diff_F_estimate} we use~\eqref{e.diff_F_2} and write for $t\in[0,T]$, $\epsilon \in \left(0, \nicefrac{(t+1)}{2} \right)$, $r \in [(t-\varepsilon) \lor 0, \, (t+\varepsilon) \land T)]$, and $s \in [t-2\epsilon, r)$:

\[
 \left| \dr F(r,s) \right| \le  |u| \int_{\R_0}x  \left( e^{iux f(r,s)} -1\right) (1+ g^*(x,s)) \ \nu(\dd x) \left| \dr f(r,s)\right| 
 \le \tilde C |s|^{-\beta}|r-s|^{-\gamma}
\]
with $\beta$ and $\gamma$ given by \eqref{e.decay_of_dtf} and $\tilde C$ defined as follows:

\[
 \tilde C:=  C_0 u^2 \sup_{y \in \R_0} \left| \frac{e^{iy}-1}{y} \right| \sup_{u\in [(t-\varepsilon) \lor 0, \, t+\varepsilon]}\sup_{v \in [t-2\epsilon,t+\epsilon]} \left| f(u,v)  \right| \int_{\R_0}x^2   (1+ g^*(x,s)) \ \nu(\dd x) < \infty,
\]
where the finiteness results from Definition~\ref{d.1}\eqref{continuity_of_f}.

Therefore, the assumptions of Lemma~\ref{l.continuity_of_integral}\eqref{F_i_conditions} are satisfied and thus we conclude by using the expression for the derivative obtained in \eqref{e.diff_F_2} that

\[
 \begin{aligned}
 \ddt &\int_{-\infty}^t \int_{\R_0} \left( e^{iuxf(t,s)} -1 -iux f(t,s) \right)(1+ g^*(x,s))  \ \nu (\dd x) \  \dd s  \\
&=\int_{\R_0} \left( e^{iuxf(t,t)} -1 -iux f(t,t) \right)(1+ g^*(x,t)) \  \nu
(\dd x) 
\\
& \quad +\int_{-\infty}^t \int_{\R_0} \left(iux  \dt f(t,s) \left( e^{iuxf(t,s)} -1
\right)(1+ g^*(x,s)) \right) \  \nu (\dd x)  \ \dd s 
 \end{aligned}
\]

and this derivative is continuous on $[0,T]$.

The continuous differentiability of the third summand in the exponential of \eqref{S(exp(iuM(t)))(g)} can be proven using Lemma~\ref{l.continuity_of_integral} by choosing

\[
 F(t,s) := f(t,s)^2
\]

in a similar but much easier way. The details are therefore omitted.
\end{pf}

The following result generalises the first part of Proposition~4.2 in~\cite{LS_ddf}.

\begin{proposition}\label{schwartzfunktion}
Suppose that $L$ has a nontrivial Gaussian part (i.e. $\sigma>0$) and that the moment condition

\[
 \int_{\R_0} x^n \ \nu(\dd x) < \infty \quad \text{for all } n \ge 2
\]

is fulfilled. Then for every $t \in (0,T]$ the mapping
$u \mapsto \E^{\qg}\left(e^{iuM(t)}\right)$ is a Schwartz function on $\R$.
\end{proposition}

\begin{pf}
The proof is divided into two parts. The first part deals with the derivative of the map
$u \mapsto \E^{\qg}\left(e^{iuM(t)}\right)$, which is then used in the second part 
to prove the assertion.

\underline{Part I}
First we show by induction that the above mapping is smooth with $j$-th derivative, $j \in
\N \cup \{0\}$, given by
\begin{equation}\label{e.j-derivative}
 \frac{\dd^j}{\dd u^j} \E^{\qg}\left(e^{iuM(t)}\right) = \E^{\qg}\left(i^j M(t)^j e^{iuM(t)}\right).
\end{equation}
For $j=0$ the assertion is trivial. Now let the statement hold for some $k \in
\N$. For the purpose of interchanging differentiation and integration we
consider
\[
 \begin{aligned}
  \sup_{h \in \R_0} \left| i^k M(t)^k \frac{e^{i(u+h)M(t)}-e^{iuM(t)}}{h}
\right| &= \left|M(t)^k\right| \sup_{h \in \R_0} \left| M(t)
\frac{e^{ihM(t)}-1}{hM(t)} \right| \cdot \left|  e^{iuM(t)}\right| \\ 
&\leq \left|M(t)^{k+1}\right| \sup_{x \in \R_0} \left| \frac{e^{ix}-1}{x} \right|.
 \end{aligned}
\]
Since the last supremum is finite, the term on the right-hand side is bounded by
$D|M(t)^{k+1}|$ for some $D> 0$. In the light of
\eqref{estimate_changed_measure} this yields
\[
 \E^{\qg}\left( \sup_{h \in \R_0} \left| i^k M(t)^k
\frac{e^{i(u+h)M(t)}-e^{iuM(t)}}{h} \right| \right) \leq \E^{\qg}\left(D\left|M(t)^{k+1}\right|\right)
\leq e_g \, D \, \E\left(\left|M(t)^{2(k+1)}\right|\right)^{\nicefrac{1}{2}},
\]
which is finite thanks to Lemma~\ref{l.moments_M}, where $e_g$ is a constant depending only on
$g$. Therefore, we can apply the DCT in order to interchange differentiation and
integration and obtain 
\[
\begin{aligned}
 \frac{\dd^{k+1}}{\dd u^{k+1}}\E^{\qg}\left(e^{iuM(t)}\right) &=  \frac{\dd}{\dd u}\E^{\qg}\left(
i^k M(t)^k e^{iuM(t)}\right) = \lim_{h \to 0} \E^{\qg}\left( i^k M(t)^k
\frac{e^{i(u+h)M(t)}-e^{iuM(t)}}{h}  \right) \\
&=  \E^{\qg}\left( i^{k+1} M(t)^{k+1} e^{iuM(t)}\right), 
\end{aligned}
\]
where the first equality follows from the induction hypothesis. This proves
(\ref{e.j-derivative}) and hence finishes the first part of the proof. 

\underline{Part II}
It remains to show that for all $m, n \in \N \cup \{0\}$ the expression 
\begin{equation}\label{derivative_schwartz_function}
  \left| u^n \frac{\dd^{m}}{\dd u^{m}}\E^{\qg}\left(e^{iuM(t)}\right) \right|
\end{equation}

is bounded in $u$. In view of Proposition~\ref{eiumt} we start by writing
\begin{equation}\label{e.decomp_char_fct}
 \E^{\qg}\left(e^{iuM(t)}\right) = \exp\left(-\frac{\sigma^2}{2} \int_{-\infty}^t f(t,s)^2
 \ \dd s \cdot u^2\right) \cdot R_{g,t}(u)
\end{equation}
with $R_{g,t}(u)$ given by
\[
 R_{g,t}(u) = \exp \left( iu \int_{-\infty}^t f(t,s) (\sigma^2  g(0,s))  \ \dd s
\right) \cdot \E^{\qg}\left(e^{iuM_j(t)}\right),
\]
where the process $M_j$ is constructed analogously to $M$ by using the characteristic triple $(\gamma, 0, \nu)$ instead of $(\gamma, \sigma, \nu)$. Applying the arguments of Part I to the process $(M_j(t))_{t\in\R}$  we infer that $R_{g,t}$ has bounded derivatives of every order. Since $\frac{\sigma^2}{2} \int_{-\infty}^t f(t,s)^2  \ \dd s > 0$ according to Definition~\ref{d.1}\eqref{no_zero_set}, the mapping 
\[
u
\mapsto \exp\left(-\frac{\sigma^2}{2} \int_{-\infty}^t f(t,s)^2 \  \dd s \cdot
u^2\right)
\]

is a Schwartz function and thus \eqref{derivative_schwartz_function} is bounded in $u$,
which completes the proof.
\end{pf}

\begin{rem}\label{derivative_schwartzfunktion}
Let $\sigma > 0$, $t\in (0, T)$ and $0< \delta < t \land (T-t)$. Note that

\[
\sup_{s\in[t-\delta,t+\delta]} \left|\exp\left(-\frac{\sigma^2}{2} \int_{-\infty}^s f(s,r)^2
 \ \dd r \cdot u^2\right) \right| = \exp\left(-\frac{\sigma^2}{2} \int_{-\infty}^{s_0} f(s_0,r)^2
 \ \dd r \cdot u^2\right) \le \exp\left(-c u^2\right)
\]

for some $s_0 \in [t-\delta,t+\delta]$ and $c>0$, cf. Definition~\ref{d.1}\eqref{no_zero_set}. Furthermore, it follows from \eqref{estimate_changed_measure} that

\begin{equation}\label{e.eqg_estimate}
\sup_{u\in\R}\sup_{r\in[t-\delta, t+\delta]} \left|\mathbb
E^{\qg}\left(e^{iuM_j(r)}\right)\right| \le \sup_{u\in\R}\sup_{r\in[t-\delta, t+\delta]} \mathbb
E\left(\left|e^{iuM_j(r) }\exp^\diamond(I_1(g))\right| \right)\le e_g .
\end{equation}

Thus, we see in view of \eqref{e.decomp_char_fct} that the function

\[
u \mapsto \sup_{s\in[t-\delta,t+\delta]} \E^{\qg}\left(e^{iuM(s)}\right)
\]

is an element of $\mathscr L^2(\dd u)$ in the situation that $L$ has a nontrivial Gaussian part.
\end{rem}

Our starting point in the proof of Theorem~\ref{main_result_L2} is the equation

\[
 S(G(M(T)))(g) = G(0)+ \int_0^T \ddt S(G(M(t)))(g) \  \dd t,
\]

which follows from the fundamental theorem of calculus if the function $t \mapsto S(G(M(t)))(g)$ is continuously differentiable. Therefore, we first prove the existence of $\ddt S(G(M(t)))(g)$. 

\begin{proposition}\label{e.derivative_SGM}
Let $f \in \mathcal K$, $g \in \Xi$ as well $G\in C^2(\R)$ such that one of the following assumptions holds:

\begin{enumerate}[a)]
 \item $\sigma>0$ and $G$ has compact support, 
 \item $G,G',G''\in\mathcal A(\R)$.
  \end{enumerate} 
  
 Then $S(G(M(\cdot)))(g)$ is continuously differentiable on $[0,T]$ with derivative
 
\[
 \ddt S(G(M(t)))(g)= \frac{1}{\sqrt{2\pi}}\int_\R \left( \F G\right)(u)\dt \mathbb
E^{\qg}\left(e^{iuM(t)}\right)\ \dd u.
\]
\end{proposition}

\begin{pf}
Let $t\in (0, T)$ and $0< \delta < t \land (T-t)$. The cases $t=0$ and $t=T$ with the one-sided limits can be handled analogously. Recall that the Fourier transform of a function $f \in \L^1(\R)$ is denoted by $\F f$. By means of the Fourier inversion theorem we deduce that 

\begin{equation}\label{S_transform_via_pairing}
  S(G(M(t)))(g) = \E^{\qg}(G(M(t))) =   \frac{1}{\sqrt{2 \pi}} \int_{\R} (\F
G)(u) \E^{\qg}\left( e^{iuM(t)} \right)  \ \dd u
\end{equation}

holds under condition $b)$. For a function $G$ fulfilling condition $a)$ we use a standard approximation via convolution with a $C^\infty$-function with compact support and deduce by means of the DCT and Remark~\ref{derivative_schwartzfunktion} that~\eqref{S_transform_via_pairing} also holds true in that case.  Hence, in both cases we have

\begin{equation}\label{e.diff_quotient}
\begin{aligned}
 &\sqrt{2\pi}\frac{S(G(M(t+h)))(g)-S(G(M(t)))(g)}{h} 
\\[0.5ex]
&= \int_\R \left( \F G\right)(u)  \frac{\mathbb
E^{\qg}\left(e^{iuM(t+h)}\right)-\mathbb
E^{\qg}\left(e^{iuM(t)}\right)}{h}\ \dd u.
\end{aligned}
\end{equation}

Moreover, resorting to  Lemma~\ref{derivative_S(exp(iuM(t)))} we infer from  the mean value theorem for complex valued functions that 

\[
\begin{aligned}
& \sup_{h\in[-\delta,\delta]\setminus\{0\}}\left|\left( \F  G\right)(u)  \frac{\mathbb
E^{\qg}\left(e^{iuM(t+h)}\right)-\mathbb
E^{\qg}\left(e^{iuM(t)}\right)}{h}\right|
\\[0.5ex]
&\le \left|\left( \F  G\right)(u)\right|\cdot \sup_{r\in[t-\delta, t+\delta]} \left|\dr \mathbb
E^{\qg}\left(e^{iuM(r)}\right)\right|
\\
&=\left|\left( \F  G\right)(u)\right|\cdot  \sup_{r\in[t-\delta, t+\delta]} \left|\mathbb
E^{\qg}\left(e^{iuM(r)}\right)\left({\rm I}(r)+{\rm II}(r)+{\rm III}(r)+{\rm IV}(r)\right)\right|,
\end{aligned}
\]

where the terms I$(r)$, II$(r)$, III$(r)$ and IV$(r)$ are given by 

\[
 \begin{aligned}  
\text{I}(r)&= iu \frac{\dd }{\dd r}  S(M(r))(g),
\\[0.5ex]
\text{II}(r)&= -\frac{\sigma^2u^2}{2}\left(f(r,r)^2+ 2\int_{-\infty}^r \dr f(r,s)\cdot f(r,s) \ 
\dd s\right) ,
\\[0.5ex]
\text{III}(r)&= \int_{\R_0} \left( e^{iuxf(r,r)} -1 -iux f(r,r) \right)(1+ g^*(x,r)) \  \nu
(\dd x) 
\\[0.5ex] 
\text{and} \quad & \\
\text{IV}(r)&= \int_{-\infty}^r \int_{\R_0} \left(iux  \dr f(r,s) \left( e^{iuxf(r,s)} -1
\right)(1+ g^*(x,s)) \right) \  \nu (\dd x)  \ \dd s. 
 \end{aligned}
\]

Hence, 

\begin{align}\label{e.estimate_diff_quotient}
 & \sup_{h\in[-\delta,\delta]\setminus\{0\}} \left|\left( \F  G\right)(u)  \frac{\mathbb
E^{\qg}\left(e^{iuM(t+h)}\right)-\mathbb
E^{\qg}\left(e^{iuM(t)}\right)}{h}\right|\notag
\\
\le  &\sup_{r\in[t-\delta, t+\delta]} \left|\mathbb
E^{\qg}\left(e^{iuM(r)}\right)\right|\notag
\\& \times \sup_{r\in[t-\delta, t+\delta]}\bigg(\left|(\mathcal F G')(u)\right|\cdot\left|\frac{\dd }{\dd r}  S(M(r))(g)\right|\notag
\\
&\hspace{2.75cm} +\left|(\mathcal F G'')(u)\right|\cdot\bigg[
-\frac{\sigma^2}{2}\left(f(r,r)^2+ 2\int_{-\infty}^r \dr f(r,s)\cdot f(r,s) \ 
\dd s\right)
\\
&\hspace{3.25cm}+\sup_{y\in\R_0}\left|\frac{e^{iy}-1-iy}{y^2}\right|\int_{\R_0} x^2f(r,r)^2\left|1+ g^*(x,r)\right| \  \nu
(\dd x) \notag
\\
&\hspace{3.25cm}+\sup_{y\in\R_0}\left|\frac{e^{iy}-1}{y}\right|\int_{-\infty}^r \int_{\R_0} \left(x^2  \left|\dr f(r,s)\right|\cdot\left|f(r,s)\right| \left|1+ g^*(x,s)\right| \right) \  \nu (\dd x)  \ \dd s\bigg]\bigg)\notag
\\
=& \sup_{r\in[t-\delta, t+\delta]} \left|\mathbb
E^{\qg}\left(e^{iuM(r)}\right)\right|\cdot\Biggl(\left|(\mathcal F G')(u)\right|\cdot\sup_{r\in[t-\delta, t+\delta]}\left|\frac{\dd }{\dd r}  S(M(r))(g)\right| \notag \\
   &\hspace{8cm}+ \left|(\mathcal F G'')(u)\right|\cdot\sup_{r\in[t-\delta, t+\delta]}\mathfrak E(r)\Biggr),\notag
\end{align}  

where $\mathfrak E(r)$ denotes the expression in the square brackets. We now prove that under assumption a) or b) the left-hand side of~\eqref{e.estimate_diff_quotient} is in  $\mathscr L^1(\dd u)$:

\begin{enumerate}[a)]
 \item If $\sigma>0$ and $G$ has compact support, then in particular $G',G'' \in \mathscr L^1 \cap \mathscr L^2 (\dd u)$ and $ \mathcal FG', \mathcal FG''\in\mathscr L^2(\dd u)$. Since $r\mapsto\mathfrak E(r)$ is continuous (cf. Lemma~\ref{derivative_S(exp(iuM(t)))}), we deduce that 
 \begin{equation}\label{e.finite_E}
 \sup_{r\in[t-\delta, t+\delta]}\mathfrak E(r)<\infty.
 \end{equation}
 According to Lemma~\ref{derivative_S(M(t))} also $r\mapsto \frac{\dd}{\dd r}  S(M(r))(g)$ is continuous and thus  
 \begin{equation}\label{e.finite_derSM}
 \sup_{r\in[t-\delta, t+\delta]}\left|\frac{\dd }{\dd r}  S(M(r))(g)\right|<\infty.
 \end{equation}
 In view of Remark~\ref{derivative_schwartzfunktion} we thus infer that both factors on the right-hand side of (\ref{e.estimate_diff_quotient}) are elements of $\mathscr L^2(\dd u)$. By means of  the Cauchy-Schwarz inequality it follows that
 \begin{equation}\label{e.L1}
 u\mapsto \sup_{h\in[-\delta,\delta]\setminus\{0\}} \left|\left( \F  G\right)(u)  \frac{\mathbb
E^{\qg}\left(e^{iuM(t+h)}\right)-\mathbb
E^{\qg}\left(e^{iuM(t)}\right)}{h}\right|\in\mathscr L^1(\dd u).
\end{equation}
 
 \item Note that if $G',G''\in\mathcal A(\R)$, then $\mathcal FG',\mathcal FG''\in\mathscr L^1(\dd u)$. 
Therefore, in the light of \eqref{e.eqg_estimate}, (\ref{e.finite_E}) and (\ref{e.finite_derSM}) we conclude that (\ref{e.L1}) holds.
\end{enumerate}

In both cases we are now able to apply the DCT to (\ref{e.diff_quotient}) and thus we infer that $t \mapsto S(G(M(t)))(g) $ is differentiable and

\[
\begin{aligned}
\ddt S(G(M(t)))(g) &= \lim_{h\to0}\frac{S(G(M(t+h)))(g)-S(G(M(t)))(g)}{h} 
\\[0.5ex]
&= \frac{1}{\sqrt{2\pi}} \int_\R \left( \F G\right)(u)  \lim_{h\to0}\frac{\mathbb
E^{\qg}\left(e^{iuM(t+h)}\right)-\mathbb
E^{\qg}\left(e^{iuM(t)}\right)}{h}\ \dd u
\\
&= \frac{1}{\sqrt{2\pi}} \int_\R \left( \F G\right)(u)\dt \mathbb
E^{\qg}\left(e^{iuM(t)}\right)\ \dd u.
\end{aligned}  
\]

The continuity of the derivative can be proven by using Lemma~\ref{derivative_S(exp(iuM(t)))} and the estimates in \eqref{e.estimate_diff_quotient}.
\end{pf}

We are now in a position to give the proof of our main result. 

{\bf{Proof of Theorem~\ref{main_result_L2}}}
The proof is divided into two parts. In the first part we show that the desired formula holds if 

\begin{enumerate}[a)]
\item $\sigma>0$ and $G$ has compact support,
or
\item $G,G',G''\in\mathcal A(\R)$.
\end{enumerate}

In the second part we are concerned with proving  Theorem~\ref{main_result_L2} (i), i.e.  we deal with the situation that $\sigma>0$ and  $G, G'$ and $G''$ are of polynomial growth.  To this end, we approximate such functions $G$ by functions that satisfy a).

\underline{Part I}
Let $G$ satisfy a) or b) above. By means of  Proposition~\ref{e.derivative_SGM} and Lemma~\ref{derivative_S(exp(iuM(t)))} we obtain

\begin{align}\label{e.key.0}
 S(G(M(T)))(g)-S(G(M(0)))(g) &= \int_0^T \ddt S(G(M(t)))(g) \  \dd t
\\[0.5ex]
 & = \frac{1}{\sqrt{2 \pi}} \int_0^T  \left({\rm I}(t)+{\rm II}(t)+{\rm III}(t)+{\rm IV}(t)\right) \ \dd t \notag
\end{align}

with the terms $\text{I}(t)$, $\text{II}(t)$, $\text{III}(t)$ and $\text{IV}(t)$ given by

\[
 \begin{aligned}  
\text{I}(t)&= \int_{\R} (\F G)(u) \E^{\qg}\left( e^{iuM(t)} \right)
\cdot  iu  \frac{\dd }{\dd t}  S(M(t))(g) \  \dd u,
\\[0.5ex]
\text{II}(t)&= \int_{\R} (\F G)(u)   \frac{-\sigma^2u^2}{2}
\E^{\qg}\left( e^{iuM(t)} \right)\left(f(t,t)^2 +
2\int_{-\infty}^{t} \dt f(t,s) f(t,s)  \ \dd s \right)  \dd u,
\\[0.5ex]
\text{III}(t)&= \int_{\R} (\F G)(u)  \E^{\qg}\left( e^{iuM(t)} \right)
\int_{\R_0} \left( e^{iuxf(t,t)} -1 -iux f(t,t) \right) 
(1+ g^*(x,t))\ \nu (\dd x)\ \dd u
\\[0.5ex] 
\text{and} \quad & \\
\text{IV}(t)&= \int_{\R} (\F G)(u) \E^{\qg}\left( e^{iuM(t)} \right) 
 \cdot \int_{-\infty}^{t} \int_{\R_0} \left(iux  \dt f(t,s)\left( e^{iuxf(t,s)} -1  \right) \right) 
\\[0.5ex] 
& \hspace*{10cm}\cdot(1+ g^*(x,s)) \ \nu (\dd x)\ \dd s \ \dd u. 
 \end{aligned}
\]

We exemplarily give the argument for ${\rm III}(t)$. For this term we have

\begin{align*}
{\rm III}(t) &=  \int_{\R_0} \int_{\R} (\F G)(u)
\left( e^{iuxf(t,t)} -1 -iux f(t,t) \right)(1+ g^*(x,t)) \cdot \E^{\qg}\left( e^{iuM(t)} \right) \ \dd u\ \nu (\dd x),
\end{align*}
 where we applied Fubini's theorem, which is justified by the estimate

\[
 \left| e^{iuxf(t,t)} -1 -iux f(t,t)  \right| \le D u^2 x^2 f(t,t)^2
\]

for some constant $D >0$, and every $u \in \R$ and $x \in \R_0$ and the assumed integrability properties on~$\mathcal F G''(u)$. By using standard manipulations of the Fourier transform we derive

\begin{align*}
\text{III}(t) &= \int_{\R_0}\int_\R \F\left(G(\cdot+ xf(t,t))-  G(\cdot) -xf(t,t)
G'(\cdot)\right)(u)
\\[0.5ex]
& \qquad \qquad \cdot \E^{\qg}\left( e^{iuM(t)} \right) (1+ g^*(x,t)) \ \dd
u \ \nu (\dd x).
\end{align*}

Consequently, by applying (\ref{S_transform_via_pairing}) we obtain

\begin{align*}
& {\rm{III}}(t)= \int_{\R_0}\sqrt{2\pi} S\left( G(M(t)+ xf(t,t))-  G(M(t)) -xf(t,t) G'(M(t)) \right)(g) \
(1+g^*(x,t))  \ \nu(\dd x) .
\end{align*}

For the terms corresponding to I$(t)$, II$(t)$ and IV$(t)$ similar techniques
apply and result in 

\[
\begin{aligned}
  \ \text{I}(t)   &=  \sqrt{2\pi}
S(G'(M(t)))(g) \ddt S(M(t))(g)
\end{aligned}
\]

and

\[
\begin{aligned}
 \text{II}(t) &=  \frac{\sigma^2}{2}  \left(f(t,t)^2 + 2\int_{-\infty}^t \dt f(t,s)\cdot
f(t,s) \  \dd s \right) \int_{\R}  (\F G'')(u)   \E^{\qg}\left( e^{iuM(t)} \right) \ \dd u
\\[0.5ex]
&=  \frac{\sigma^2}{2}  \left(f(t,t)^2 + 2\int_{-\infty}^t \dt f(t,s)\cdot
f(t,s) \  \dd s \right) \sqrt{2\pi}  S(G''(M(t)))(g)
\end{aligned}
\]

as well as

\[
\begin{aligned}
 &  \ \text{IV}(t) = \int_{-\infty}^t \int_{\R_0} x \dt f(t,s) \sqrt{2\pi} S(G'(M(t)+ xf(t,s))- G'(M(t)))(g)
\cdot (1+ g^*(x,s))\ \nu (\dd x) \  \dd s.
\end{aligned}
\]

Therefore, plugging the above expressions into \eqref{e.key.0} we obtain

\begin{align*}
& S(G(M(T)))(g)-S(G(M(0)))(g)\notag 
\\[0.5ex]
&= \int_0^T S(G'(M(t)))(g) \ddt S(M(t))(g) \ \dd t
\\[0.5ex]
&\quad\qquad +  \int_0^T \frac{\sigma^2}{2}  \left(f(t,t)^2 + 2 \int_{-\infty}^t \dt
f(t,s)\cdot f(t,s) \  \dd s \right)  S(G''(M(t)))(g) \ \dd t\notag
\\[0.5ex]
&\quad\qquad + \int_0^T \int_{\R_0}  S\Bigl( G(M(t)+ xf(t,t))- G(M(t))
-xf(t,t)G'(M(t))\Bigr)(g) \cdot (1+ g^*(x,t))\ \nu (\dd x) \ \dd t\notag
\\[0.5ex]
&\quad\qquad +\int_0^T \int_{-\infty}^t \int_{\R_0} x \dt f(t,s) S(G'(M(t)+ xf(t,s))-
G'(M(t)))(g) \cdot (1+ g^*(x,s))\ \nu (\dd x) \  \dd s \ \dd t\notag
\\[0.5ex]
& =:{\rm  {I^*}}+{\rm { II^*}}+{\rm{ III^*}}+{\rm{ IV^*}}.\notag
\end{align*}

An application of Lemma~\ref{derivative_S(M(t))}  and Fubini's theorem yields for ${\rm I}^*$ the following equality

\[
\begin{aligned}
 &\int_0^T S(G'(M(t)))(g) \ddt S(M(t))(g) \ \dd t \\ 
 &= \int_{-\infty}^T \int_{0 \lor s}^T S(G'(M(t)))(g)  \sigma   \dt f(t,s) g(0,s) \ \dd t \ \dd s  \\
 &\quad+ \int_0^T S(G'(M(t)))(g)  \sigma f(t,t) g(0,t) \ \dd t \\ 
 &\quad+  \int_0^T S(G'(M(t)))(g) f(t,t) \int_{\R_0}x g^*(x,t) \ \nu (\dd x) \ \dd t \\
 &\quad+ \int_0^T S(G'(M(t)))(g)  \int_{-\infty}^t \int_{\R_0} \dt f(t,s) x g^*(x,s)\ \nu(\dd x) \ \dd s \ \dd t \\
 &  =:\text{I}^{*}_a+ \text{I}^{*}_b+ \text{I}^{*}_c+ \text{I}^{*}_d.
\end{aligned}
\]

Resorting to Remark~\ref{special_cases} we obtain

\[
 \text{I}^{*}_b = S\left(\sigma \int_0^T G'(M(t-)) f(t,t) \ W(\dd t) \right)(g)
\]

as well as

\[
 \text{I}^{*}_c = S\left(\int_0^T \int_{\R_0} G'(M(t-)) f(t,t) x \ \tilde{N}(\dd x, \dd t) \right)(g),
\]

where we used that both integrals exist separately and reduce to classical stochastic integrals, because $t$ is $\qg$-a.s. not a jump time and therefore $S(G'(M(t-)))(g)=S(G'(M(t)))(g)$ and because of the predictability of the integrands.
In particular, we infer 

\[
\text{I}^{*}_b+\text{I}^{*}_c = 
S\left(\int_0^T G'(M(t-)) f(t,t) \ L(\dd t) \right)(g).
\]

In view of the boundedness of $G''$ we can use Fubini's theorem to write the term $\text{II}^*$ as

\[
\begin{aligned}
   {\text{II}}^*   =   S\left( \frac{\sigma^2}{2} \int_0^T     \left(f(t,t)^2 + 2 \int_{-\infty}^t
\dt f(t,s)\cdot f(t,s)  \ \dd s \right)  G''(M(t))   \   \dd t \right)(g). 
\end{aligned}
\]

For the third term we obtain 

 \begin{align*}
   {\text{III}}^* &= S \left(\int_0^T      \int_{\R_0}   \left( G(M(t-)+ xf(t,t))- G(M(t-))
-xf(t,t)G'(M(t-)) \right)\ N^\diamond (\dd x, \dd t) \right) (g) \notag
\\
&= S \left( \sum_{0<t \leq T} G(M(t))- G(M(t-))- G'(M(t-))\Delta M(t)\right) (g).
\end{align*}

Note that the predictability of the integrand implies that the Hitsuda-Skorokhod
integral is an ordinary integral with respect to the jump measure and hence the second equality above holds by means
of Lemma~\ref{l.moments_M}.

Using the boundedness of $G'$ and $G''$, the mean value theorem and Lemma~\ref{l.continuity_of_integral}(i), with $F(t,s) = \left|f(t,s ) \ds f(t,s) \right|$, in order to justify the application of Fubini's theorem we deduce that

\[
\begin{aligned}
  {\text{IV}}^* &=\int_0^T \int_{-\infty}^t \int_{\R_0} x \dt f(t,s) S(G'(M(t)+ xf(t,s))- G'(M(t)))(g) \cdot (1+ g^*(x,s))\ \nu (\dd x) \  \dd s \ \dd t  \\
  &= \int_{-\infty}^T \int_{\R_0} \int_{0 \lor s}^T x \dt f(t,s) S \left(
G'(M(t)+ xf(t,s) \right)(g) \cdot g^*(x,s) \ \dd t \ \nu(\dd x) \ \dd s   \\
  &\quad-\int_0^T    \int_{-\infty}^t \int_{\R_0} x \dt f(t,s) S(G'(M(t))(g) \cdot g^*(x,s)  \ \nu(\dd x) \ \dd s  \ \dd t
   \\
&\quad +S \left(\int_{-\infty}^T \int_{\R_0} \int_{0 \lor s}^T x \dt f(t,s) \left(
G'(M(t)+ xf(t,s))- G'(M(t)) \right)\ \dd t \  \nu (\dd x) \ \dd s  \right)(g) \\
&=: \text{IV}^{*}_a - \text{I}^{*}_d+ \text{IV}^{*}_c.
\end{aligned}
\]

Therefore we obtain in view of \eqref{e.key.0} that

\begin{equation}\label{e.key.0.b}
 \begin{aligned}
  \text{I}^{*}_a&+ \text{IV}^{*}_a \\
  &=S(G(M(T))-G(0))(g) \\
&\qquad  -S\left(  \frac{\sigma^2}{2} \int_0^T  \left(f(t,t)^2 + 2 \int_{-\infty}^t
\dt f(t,s)\cdot f(t,s)  \ \dd s \right)  G''(M(t))   \   \dd t \right)(g)
\\
&\qquad -S \left( \sum_{0<t \leq T} G(M(t))- G(M(t-))- G'(M(t-))\Delta M(t)\right) (g)
\\
&\qquad -S \left( \int_{-\infty}^T \int_{\R_0} \int_{0 \lor s}^T x \dt f(t,s) \left(
G'(M(t)+ xf(t,s))- G'(M(t)) \right)\ \dd t \   \nu (\dd x) \ \dd s  
\right) (g) 
\\ 
&\qquad
-S\left(\int_0^T G'(M(t-)) f(t,t) \ L(\dd t) \right)(g),
 \end{aligned}
\end{equation}

where the existence in $\L^2(\mathbb P)$ of the arguments of the $S$-transforms on the right-hand side follows from the boundedness of $G$, $G'$ and $G''$, using Lemma~\ref{l.moments_M} and the $\mathscr L^4(\mathbb P)$-assumptions on $L$ for the sum over the jumps of $M$.

By linearity of the $S$-transform I$^{*}_a$ $+$ IV$^{*}_a$ equals the $S$-transform of some $\Phi\in\mathscr L^2(\mathbb P)$ and can thus, by Definition~\ref{Lambda-diamond}, be written as

\[
 S\left(\int_{-\infty}^T \int_{\R} \int_{0 \lor s}^T  \dt f(t,s)  \left(
G'(M(t)+ xf(t,s) \right)  \ \dd t \ \Lambda^\diamond(\dd x,\dd s) \right)(g).
\]

Hence, reordering the terms in  \eqref{e.key.0.b} and resorting to the injectivity property of the $S$-transform (see Lemma~\ref{injectivity_s-transform}) this results in the change of variable formula
\eqref{ito_formula}.

\underline{Part II} We now consider the case that $\sigma>0$ and  $G, G'$ and $G''$ are of polynomial growth of degree ${q}$, that is we have

\[
 |G^{(l)}(x)| \le C_{pol}(1+|x|^{q}) \quad \text{ for every $x \in \R$ and $l=0,1,2$}.
\]

Our approach is to approximate such $G$ by functions fulfilling condition a) above as we know from Part~I of this proof that the generalised It\=o formula holds true for these functions. To this end, let $(G_n)_{n\in\N}$ be a sequence of functions defined by 

\[
 G_n(x)= \begin{cases}
          G(x), \quad &|x| \le n, \\
          G(x) \varphi(|x|-n) , \quad &n< |x| < n+1, \\
          0, \quad &|x| \ge n+1,
         \end{cases}
\]

where the function $\varphi:(0,1) \to \R$ is given by

\[
 \varphi(y) = \exp \left(-\left( \frac{y}{1-y} \right)^3 \right).
\]

Then $G_n$ satisfies the following conditions:

\begin{enumerate}
 \item $G_n \in  C^2(\R)$,
 \item $G_n(x)= G(x), \ x \in [-n,n]$,
 \item $G_n$ has compact support and
 \item $|G_n(x)|+ |G'_n(x)|+ |G''_n(x)| \le \tilde{C}(1+|x|^{q})$ for all $x \in \R$,
\end{enumerate}

where the constant $\tilde{C}$ is given by 

\[
 \tilde{C}:= 3C_{pol} \max \left\{ \sup_{y \in (0,1)} \varphi(y) , \sup_{y \in (0,1)} \varphi'(y), \sup_{y \in (0,1)} \varphi''(y)\right\}.
\]

Below we denote by (\ref{ito_formula})$_n$ and (\ref{e.key.0.b})$_n$ the formulas~(\ref{ito_formula}) and (\ref{e.key.0.b}) with $G$ replaced by $G_n$. In Part~I we showed that the generalised It\=o formula~(\ref{ito_formula})$_n$ holds for any $n\in\N$.
It remains to show that this formula also holds true in the limit as $n \to \infty$. For this purpose we need to interchange the limits and the integrals on the right-hand side of~(\ref{ito_formula})$_n$.

Let us first deal with the penultimate term  on the right-hand side of (\ref{ito_formula})$_n$. For this purpose we use the mean value theorem to find some $z_0\in\bigl(M(t)\land (M(t)+xf(t,s)),\,M(t)\lor (M(t)+xf(t,s))\bigr)$ such that

\begin{equation}\label{e.mean_value_polynomial}
 \begin{aligned}
|G_n'(M(t)+xf(t,s))-G_n'(M(t))| &= |G_n''(z_0)|\cdot|xf(t,s)|\\ 
  &\le \tilde C (1+|z_0|^{q})\cdot|xf(t,s)| \\
  & \le \tilde C(1+ |M(t)|^{q} + |M(t)+ x f(t,s)|^{q})\cdot|xf(t,s)|.
\end{aligned}
\end{equation}

In order to apply the DCT we introduce the abbreviation

\[
 D^1_n(x,t,s) := G'(M(t)+ xf(t,s)) -G'(M(t)) - \left(G_n'(M(t)+ xf(t,s)) -G_n'(M(t))\right).
\]

for $x \in \R_0$ and $s,t \in \R$. By using \eqref{e.mean_value_polynomial} and Jensen's inequality we get the following chain of estimates

\begin{equation}\label{e.2}
 \begin{aligned}
    &\E \left(  \left(\int_{0}^T \int_{-\infty}^t \int_{\R_0}\sup_{n \in \N} \left| D^1_n(x,t,s) x \dt f(t,s) \right| \ \nu(\dd x)\  \dd s \ \dd t \right)^2 \right)  \\
    &\le \E \left(  \left(\int_{0}^T \int_{-\infty}^t \int_{\R_0} 2 \tilde C  \left( 1+ |M(t)|^{q} + |M(t)+ x f(t,s)|^{q} \right) \right. \cdot \left. \left|x^2f(t,s)  \dt f(t,s) \right| \ \nu(\dd x)\  \dd s \ \dd t \right)^2 \right) \\
    &\le  \E \left(  \left(\int_{0}^T \int_{-\infty}^t \int_{\R_0} K_1 \left( 1 + |M(t)|^{q} + |x f(t,s)|^{q} \right) \right. \cdot \left. \left|x^2f(t,s)  \dt f(t,s) \right| \ \nu(\dd x)\  \dd s \ \dd t \right)^2 \right) \\
    &\le K_1^2\E \Bigg(\bigg(   \sup_{t \in [0,T]} \left( 1 + |M(t)|^{q} \right)\int_{0}^T \int_{-\infty}^t \int_{\R_0}     \left|x^2f(t,s)  \dt f(t,s) \right| \ \nu(\dd x)\  \dd s \ \dd t \\
    & \hspace{4cm}+   \int_{0}^T \int_{-\infty}^t \int_{\R_0}   |x f(t,s)|^{q}   \cdot  \left|x^2f(t,s)  \dt f(t,s) \right| \ \nu(\dd x)\  \dd s \ \dd t  \bigg)^2 \Bigg) \\
    &\le 2K_1^2 \E   \left( \left( 1 + \sup_{t \in [0,T]}|M(t)|^{q} \right)^2 \right) \left( \int_{0}^T \int_{-\infty}^t \int_{\R_0}     x^2\left|f(t,s)  \dt f(t,s) \right| \ \nu(\dd x)\  \dd s \ \dd t \right)^2 \\
    &\hspace{4cm}+2K_1^2 \left(\int_{0}^T \int_{-\infty}^t\int_{\R_0}   |x|^{{q}+2} |f(t,s)|^{{q}+1}   \cdot  \left| \dt f(t,s) \right| \ \nu(\dd x)\  \dd s \ \dd t  \right)^2
 \end{aligned}
\end{equation}

for a suitable constant $K_1>0$. 

We now want to apply Lemma~\ref{l.continuity_of_integral}\eqref{F_i_conditions_2} to prove that the last term of \eqref{e.2} is finite. For that purpose we define for any $t\in[0,T]$, $s\le t$ and $q \ge 0$

\[
 F_{q}(t,s) :=  \int_{\R_0}   |x|^{{q}+2} |f(t,s)|^{{q}+1}   \cdot  \left| \dt f(t,s) \right| \ \nu(\dd x).
\]

Estimate~\eqref{e.decay_of_dtf} then yields

\[
\begin{aligned}
 |F_{q}(r,s)| &\le C_0  \int_{\R_0}   |x|^{{q}+2}  \ \nu(\dd x) \cdot |f(r,s)|^{{q}+1}   \cdot  \left| s \right|^{-\beta} \left| r-s \right|^{- \gamma}
\end{aligned}
\]
for every $t\in[0,T]$, $\epsilon\in(0,\nicefrac{t}{2})$, $r\in[(t-\epsilon)\lor0,\,(t+\epsilon)\land T]$, and $s< r$.
In particular, \eqref{e.F_estimate} is fulfilled with

\[
 \tilde C := C_0 \sup_{(r,v) \in [0,T]\times[0,T]}|f(r,v)|^{{q}+1} \int_{\R_0}   |x|^{{q}+2}  \ \nu(\dd x).
\]

By means of Definition~\ref{d.1}\eqref{continuity_of_f}, for all $s \in \R$ the mapping $t \mapsto F_q(t,s)$ is continuous on $(s,\infty)$. Moreover,  for any fixed $t\in[0,T]$ as well as  $C_0,\beta,\gamma, \theta$ and $\epsilon \in \left(0, \nicefrac{(t+1)}{2} \right)$ such that~\eqref{e.decay_of_dtf} and~\eqref{e.decay_of_f} hold,  we infer in view of Definition~\ref{d.1}\eqref{continuity_of_f} that

\[
 \begin{aligned}
  &\int_{-\infty}^{t-2 \epsilon} \sup_{r \in [0 \lor(t-\epsilon),\, (t+\epsilon) \land T]} \left|F_{q}(r,s) \right| \ \dd s \\
  &\le C_0   \int_{\R_0}   |x|^{{q}+2}  \ \nu(\dd x)\left(\sup_{u \in (-\infty,-1]} \sup_{r \in [0 \lor(t-\epsilon),\, (t+\epsilon) \land T]} \Bigl(|f(r,u)||u|^{\theta}\Bigr)^{{q}+1}\right) \int_{-\infty}^{-1} |s|^{-(({q}+1)\theta+\beta)} |t-s|^{-\gamma} \  \dd s \\
&\qquad +  C_0  \int_{\R_0}   |x|^{{q}+2}  \ \nu(\dd x) \cdot\sup_{u \in [-1,t]}\sup_{r\in[0 \lor(t-\epsilon),\, (t+\epsilon) \land T]}|f(r,u)|^{{q}+1} \int_{-1}^{t-2 \epsilon} |s|^{-\beta} |t-s|^{-\gamma}  \  \dd s \\
&< \infty,
 \end{aligned}
\]

and hence \eqref{e.F_integral_bounded} holds. Note that here we have used the assumptions on the moments of the L\'evy process $L$.

In view of the continuity of $f$ and $\dt f(\cdot,s)$ for Lebesgue-a.e. $s\in\R$, condition~ \eqref{F_i_conditions_2}a) in Lemma~\ref{l.continuity_of_integral} is satisfied and consequently Lemma~\ref{l.continuity_of_integral} implies that for every $q \ge 0$ the mapping

\[
 t \mapsto \int_{-\infty}^t F_{q}(t,s) \ \dd s
\]

is continuous and thus

\[
 \int_0^T  \int_{-\infty}^t F_{q}(t,s) \ \dd s \ \dd t \le T \sup_{t \in [0,T]}  \int_{-\infty}^t F_{q}(t,s) \ \dd s < \infty.
\]

Combining this with Lemma~\ref{l.moments_M} and \eqref{e.2} we see that

\begin{equation}\label{e.finiteness}
 \E \left(  \left(\int_{0}^T \int_{-\infty}^t\int_{\R_0}\sup_{n \in \N} \left| D^1_n(x,t,s) x \dt f(t,s) \right| \ \nu(\dd x)\  \dd s \ \dd t \right)^2 \right)  < \infty.
\end{equation}

Furthermore, observe that

\[
 \left|D^1_n(x,t,s)\right|\to0
\]

$\mathbb P$-a.s. for every $x \in \R_0$, $t\ge 0$ and $s \in \R$ as $n\to\infty$. In the light of \eqref{e.finiteness}  we conclude by means of Fubini's theorem and the DCT

\begin{equation}\label{e.fubini_dct}
\begin{aligned}
 \lim_{n\to\infty}& \int_{-\infty}^T \int_{\R_0} \int_{0 \lor s}^T \left( G_n'\left(M(t) + xf(t,s)\right) -G_n'(M(t)) \right) x \dt f(t,s)\ \dd t \ \nu(\dd x)\  \dd s \\
 &=\lim_{n\to\infty} \int_0^T \int_{-\infty}^t\int_{\R_0} \left( G_n'\left(M(t) + xf(t,s)\right) -G_n'(M(t)) \right) x \dt f(t,s) \ \nu(\dd x)\  \dd s \ \dd t\\
&=\int_{-\infty}^T \int_{\R_0} \int_{0 \lor s}^T  \left( G'\left(M(t) + xf(t,s)\right)-G'(M(t)) \right)x\dt f(t,s) \ \dd t \ \nu(\dd x)\  \dd s 
\end{aligned}
\end{equation}

in $\L^2(\mathbb P)$. 

To handle the third term of the right-hand side of (\ref{ito_formula})$_n$ we use Lemma~\ref{l.moments_M} and the shorthands

\[
 \delta G_n(x) := G(x)- G_n(x) \text{ and } \delta G'_n(x) := G'(x)- G'_n(x)
\]

to write

\begin{equation}\label{e.L2_difference_sum}
 \begin{aligned}
  &\E \Biggl(\Biggl(\sum_{0<t\leq T}\left[ \delta G(M(t))- \delta G(M(t-))- \delta G'(M(t-)) \Delta M(t) \right]\Biggr)^2\Biggr)\\
  &=\E \Biggl(\Biggl( \sum_{0 <t \le T} \left[\delta  G_n(M(t-)+\Delta L(t) f(t,t))-\delta G_n(M(t-)) -\delta G_n'(M(t-))\Delta L(t)f(t,t) \right]\Biggr)^2\Biggr) \\
  &=\E \left(\left(\int_0^T \int_{\R_0} \left[\delta  G_n(M(t-)+xf(t,t))-\delta G_n(M(t-)) -\delta G_n'(M(t-))xf(t,t) \right] \ N(\dd x, \dd t)\right)^2 \right)
  \end{aligned}
\end{equation} 

After introducing the abbreviation

\begin{equation}\label{d.D2}
 \begin{aligned}
 D^2_n(x,t) &:=  \delta G_n(M(t-)+xf(t,t))-\delta G_n(M(t-)) -\delta G_n'(M(t-))xf(t,t)
 \end{aligned}
\end{equation}
  
for $x \in \R_0$ and $t >0$ we continue in view of the It\=o isometry for L\'evy processes and Jensen's inequality 

\begin{equation}\label{e.difference_Dn2}
\begin{aligned}
  &\E \left(\left(\int_0^T \int_{\R_0} D^2_n(x,t) \ N(\dd x, \dd t)\right)^2 \right) \\
  &\qquad\le  2\E \left(\left( \int_0^T \int_{\R_0} D^2_n(x,t) \ \tilde N(\dd x, \dd t) \right)^2 \right) +2\E \left(\left(\int_0^T \int_{\R_0} D^2_n(x,t) \ \nu(\dd x)\ \dd t\right)^2 \right) \\
  &\qquad= 2 \int_0^T \int_{\R_0}\E \left(D^2_n(x,t)^2 \right) \ \nu(\dd x)\ \dd t +2\E \left(\left(\int_0^T \int_{\R_0} D^2_n(x,t) \ \nu(\dd x)\ \dd t\right)^2 \right) 
 \end{aligned}
\end{equation}

In order to apply the DCT we thus consider the two expressions

\[
 \text{I} := \int_0^T \int_{\R_0}\E \left(  \sup_{n \in \N} D^2_n(x,t)^2 \right)\ \nu(\dd x)\ \dd t
\]

and

\[
 \text{II} :=\E \left( \left( \int_0^T \int_{\R_0} \sup_{n \in \N} D^2_n(x,t)\ \nu(\dd x)\ \dd t \right)^2 \right).
\]

Note that by making use of the mean value theorem we can find  $z_1,z_{1,n} \in \bigl(M(t-) \land  (M(t-)+ xf(t,t)),M(t-) \lor  (M(t-)+ xf(t,t))\bigr)$ such that

\[
 \begin{aligned}
  G(M(t-)+xf(t,t))-G(M(t-)) = G'(z_1)xf(t,t)
 \end{aligned}
\]

and

\[
 \begin{aligned}
  G_n(M(t-)+xf(t,t))-G_n(M(t-)) = G_n'(z_{1,n})xf(t,t).
 \end{aligned}
\]

Plugging these expressions into \eqref{d.D2}, using the polynomial bound of $G$ and $G_n$ as well as Jensen's inequality yields

\begin{equation}\label{e.D2}
 \begin{aligned}
 \sup_{n \in \N} \left|D^2_n(x,t)\right| &= \sup_{n \in \N}\left|G'(z_1)xf(t,t)-G'(M(t-))xf(t,t) -G_n'(z_{1,n})xf(t,t) +G_n'(M(t-))xf(t,t) \right| \\
  &\le K_2|x f(t,t)| \left( 1+ |M(t-)|^{q}+ |xf(t,t)|^{q} \right)
 \end{aligned}
\end{equation}

for some constant $K_2 >0$. In view of this estimate we continue by applying Jensen's inequality

\[
 \begin{aligned}
  \text{I} &= \int_0^T \int_{\R_0}  \E \left( \sup_{n \in \N}D^2_n(x,t)^2 \right)\ \nu(\dd x)\ \dd t \\
 &\le K_2^2\int_0^T \int_{\R_0} x^2 f(t,t)^2 \E \left( \left( 1+ |M(t)|^{q}+ |xf(t,t)|^{q} \right)^2 \right)\ \nu(\dd x)\ \dd t \\
 & \le K_3 \int_0^T \int_{\R_0} x^2 f(t,t)^2 \E  \left( 1+ \sup_{r \in(0,T]}|M(r)|^{2{q}}+|xf(t,t)|^{2{q}} \right)\ \nu(\dd x)\ \dd t \\
 &< \infty
 \end{aligned}
\]

for a constant $K_3 > 0$ by using the estimate \eqref{e.D2}, where the finiteness follows again from the moment assumptions on $L$.

To handle II we make a Taylor expansion of first order of $G(M(t-)+xf(t,t))$ at $M(t-)$ to find some $z_2, z_{2,n} \in\bigl(M(t-) \land  (M(t-)+ xf(t,t)),M(t-) \lor  (M(t-)+ xf(t,t))\bigr)$ such that

\[
 G(M(t-)+xf(t,t))-G(M(t-)) -G'(M(t-))xf(t,t) = \frac{1}{2} G''(z_2) x^2 f(t,t)^2
\]

and

\[
 G_n(M(t-)+xf(t,t))-G_n(M(t-)) -G_n'(M(t-))xf(t,t) = \frac{1}{2} G_n''(z_{2,n}) x^2 f(t,t)^2.
\]

Plugging these into \eqref{d.D2}, using the polynomial bound of $G''$ and $G_n''$ as well as again Jensen's inequality results in

\[
 \begin{aligned}
 \sup_{n \in \N} \left|D^2_n(x,t)\right|   &\le K_4 x^2 f(t,t)^2 \left( 1+ |M(t-)|^{q}+ |xf(t,t)|^{q} \right)
 \end{aligned}
\]

for some constant $K _4>0$. Using this estimate leads to

\[
 \begin{aligned}
 \text{II} &\le  \E \left(  \left(K_4\int_0^T \int_{\R_0} x^2 f(t,t)^2 \E \left( 1+ |M(t)|^{q}+|xf(t,t)|^{q} \right)\ \nu(\dd x)\ \dd t \right)^2 \right) \\
 &< \infty.
 \end{aligned}
\]

Since

\[
 \left|D^2_n(x,t)\right|\to0
\]

$\mathbb P$-a.s. for every $x \in \R_0$ and $t\ge 0$ as $n\to\infty$ we see by means of the DCT that the right-hand side of \eqref{e.difference_Dn2} converges to $0$ and according to \eqref{e.L2_difference_sum} we conclude

\[
\begin{aligned}
 \lim_{n\to\infty}& \sum_{0<t\leq T} \left[ G_n(M(t))- G_n(M(t-))- G_n'(M(t-)) \Delta M(t) \right]\\
&= \sum_{0<t\leq T} \left[ G(M(t))- G(M(t-))- G'(M(t-)) \Delta M(t) \right]
\end{aligned}
\]

in $\L^2(\mathbb P)$.

Similar but simpler arguments lead to

\[
\begin{aligned}
 \lim_{n\to\infty}& \frac{\sigma^2}{2} \int_{0}^T G_n''(M(t)) \left( \ddt
\int_{-\infty}^t f(t,s)^2  \dd s \right) \dd t = \frac{\sigma^2}{2} \int_{0}^T G''(M(t)) \left( \ddt
\int_{-\infty}^t f(t,s)^2  \dd s \right) \dd t,
\end{aligned}
\]

\[
 \lim_{n\to\infty} \int_{0}^T G'_n(M(t-)) f(t,t) \ L (\dd t) =\int_{0}^T G'(M(t-)) f(t,t) \ L (\dd t)
\]
and
\[
\begin{aligned}
 \lim_{n\to\infty}& G_n(M(T))= G(M(T))
\end{aligned}
\]

in $\L^2(\mathbb P)$ as well as

\[
\begin{aligned}
 \lim_{n\to\infty}& G_n(0)= G(0).
\end{aligned}
\]

It remains to consider the $\Lambda^\diamond$-integral. Note in view of (\ref{estimate_changed_measure}) and the convergence in $\mathscr L^2(\mathbb P)$ shown above that the right-hand side of  (\ref{e.key.0.b})$_n$ converges to the right-hand side of (\ref{e.key.0.b}) as $n\to\infty$. By applying arguments as in~\eqref{e.2} and~\eqref{e.fubini_dct} we can deduce that

\[
\begin{aligned}
  & \lim_{n \to \infty} S\left(\int_{-\infty}^T \int_{\R} \int_{0 \lor s}^T \dt f(t,s)  \left(
G_n'(M(t)+ xf(t,s) \right)\ \dd t \ \Lambda^\diamond(\dd x,\dd s)  \right)(g) \\
&\qquad =S\left( \int_{-\infty}^T \int_{\R} \int_{0 \lor s}^T \dt f(t,s)  \left(
G'(M(t)+ xf(t,s) \right)\ \dd t \ \Lambda^\diamond(\dd x,\dd s) \right)(g),
\end{aligned}
\]

which in view of the linearity of the $S$-transform completes the proof of (\ref{ito_formula}). 

In order to show the variant (\ref{ito_formula_L2}), we apply the $S$-transform to the right-hand side of \eqref{ito_formula} and \eqref{ito_formula_L2}, respectively, and make use of Theorem~\ref{r.connection_m_diamond_lamnda_diamond}. After rearranging the resulting terms and using the injectivity of the $S$-transform (cf. Proposition~\ref{injectivity_s-transform}) we deduce that both right-hand sides coincide, provided all the terms exist in $\L^2(\mathbb P)$.
\hfill$\square$

\section*{Appendix}
Here we provide the proofs of Lemma~\ref{fractional_kernel_in_k} and Lemma~\ref{l.continuity_of_integral}.

{\bf Proof of Lemma~\ref{fractional_kernel_in_k}}
It is easy to see that $f_d$ satisfies \eqref{volterra}-\eqref{no_zero_set} in Definition~\ref{d.1}. Condition \eqref{behaviour_kernel} follows from the following expression containing the derivative with respect to the first argument of $f_d$, i.e.

 \[
 \dt f_d(t,s) = \frac{1}{\Gamma(d+1)}d(t-s)_+^{d-1}.
\]

If we choose  $\beta =0$ and $\gamma = 1-d$ we deduce that  equation \eqref{e.decay_of_dtf} in condition \eqref{behaviour_kernel} is satisfied. Moreover, we infer by means of the mean value theorem that 

\[
 \sup_{r \in [0 \lor(t-\epsilon),\, (t+\epsilon) \land T]} |f(r,s)| \le \sup_{r \in [0 \lor(t-\epsilon),\, (t+\epsilon) \land T]} \sup_{u \in (-s, r-s)} \frac{1}{\Gamma(d+1)} d r|u|^{d-1} \le \frac{1}{\Gamma(d+1)} d (t+ \epsilon)|s|^{d-1}.
\]

Hence, equation \eqref{e.decay_of_f} in condition \eqref{behaviour_kernel} holds with $\theta= 1-d$. That condition~\eqref{derivative_second_argument} if fulfilled follows from Example~9 in \cite{BKO_ppm}. 
\hfill$\square$

We finish this paper with the proof of Lemma~\ref{l.continuity_of_integral}.

{\bf Proof of Lemma~\ref{l.continuity_of_integral}}

(i) We start with the right continuity and therefore can assume $t\in [0,T)$. With $\varepsilon>0$ as in assumption~\eqref{F_i_conditions_2}b) we write 
 
 \[
  \begin{aligned}
   &I_F(t+h) -I_F(t) \\
   &=  \int_{t}^{t+h} F(t+h,s) \ \dd s + \int_{-\infty}^{t-2 \varepsilon} \left(F(t+h,s)- F(t,s)\right) \ \dd s + \int_{t-2 \varepsilon}^t \left(F(t+h,s)- F(t,s)\right) \ \dd s \\
   &=: \text{I}_{\rm{i}}(h)+ \text{II}_{\rm{i}}(h)+ \text{III}_{\rm{i}}(h).
  \end{aligned}
 \]
 
For the first term we use the substitution $v := \frac{s}{t+h}$, which in view of (\ref{e.F_estimate}) yields

\begin{equation}\label{e.substitution}
 \begin{aligned}
  \sup_{h \in (0,\varepsilon \land (T-t))}&\int_{t}^{t+h} \left|F(t+h,s) \right|^q \ \dd s\\
 &\le  \tilde C^q \sup_{h \in (0,\varepsilon \land (T-t))} \int_{t}^{t+h} s^{-\beta q} \left( t+h-s \right)^{-\gamma q} \ \dd s  \\
  &=  \tilde C^q \sup_{h \in (0,\varepsilon \land (T-t))} \left( (t+h)^{1-(\beta +\gamma) q} \int_{\frac{t}{t+h}}^{1} v^{-\beta q} \left(1-v \right)^{-\gamma q} \ \dd v \right)\\
  &< \infty
 \end{aligned}
\end{equation}

 for $1< q < (\beta+ \gamma)^{-1}$. Since $\mathds 1_{[t,t+h]}(s) F(t+h,s)$ converges to  $0$ pointwise for Lebesgue-a.e. $s$ as $h \downarrow 0$, an application of the de la Vall\'ee-Poussin theorem results in 

\[
  \lim_{h \downarrow 0} \text{I}_{\rm{i}}(h) = 0.
\]

For the term II$_{\rm i}(h)$ we consider

\[
 \int_{-\infty}^{t-2 \varepsilon} \sup_{h \in (0,\varepsilon \land (T-t))} \left|F(t+h,s)- F(t,s)\right| \ \dd s  
 \le 2\int_{-\infty}^{t-2 \varepsilon} \sup_{r \in [t, (t+\varepsilon)\land T]} \left|  F(r,s) \right| \ \dd s < \infty, 
\]

where the finiteness follows from \eqref{e.F_integral_bounded}. Therefore, we can apply the DCT to obtain by means of \eqref{F_i_conditions_2}a) that

\[
 \lim_{h \downarrow 0} \left(\text{II}_{\rm{i}}(h)-\int_{-\infty}^{t-2 \varepsilon}  F(t,s) \ \dd s\right)\to0.
\]

For the term III$_{\rm i}(h)$ we estimate in the case $t>0$:

\[
 \sup_{h \in (0,\varepsilon \land (T-t))} \int_{t-2 \varepsilon}^{t} \left|F(t+h,s) -F(t,s)\right|^q \ \dd s \le  2^{q} \sup_{r \in [t, (t+\varepsilon) \land T]} \int_{t-2 \varepsilon}^{t}  \left|  F(r,s) \right|^q \ \dd s < \infty,
\]

where the finiteness follows by the same considerations as in \eqref{e.substitution} with the substitution $v=\nicefrac{s}{r}$. If $t=0$, we estimate by using \eqref{e.F_estimate}:

 \[
 \begin{aligned}
  \sup_{h \in (0,\varepsilon \land T)} &\int_{-2 \varepsilon}^{0} \left|F(h,s) -F(0,s)\right|^q \ \dd s  \\
  &\le   \tilde C 2^q  \sup_{h \in (0,\varepsilon \land T)}  \left( \int_{-2 \epsilon}^0 |-s|^{- \beta q} |h-s|^{-\gamma q}  \ \dd s + \int_{-2 \epsilon}^0  |-s|^{- \beta q} |-s|^{-\gamma q}  \ \dd s \right) \\
  & \le \tilde C 2^{q+1} \int_{-2 \epsilon}^0  |-s|^{-(\beta+\gamma) q}  \ \dd s,
  \end{aligned}
 \]

which is finite for  $1< q < (\beta+ \gamma)^{-1}$. Another application of the de la Vall\'ee-Poussin theorem therefore yields the convergence

\[
  \lim_{h \downarrow 0} \left(\text{III}_{\rm{i}}(h)-\int_{t-2 \varepsilon}^t  F(t,s) \ \dd s\right)\to0,
\]

which shows the right continuity of $I_F$ at $t$. The left continuity follows by similar arguments.

(ii) We prove differentiability from above and note that differentiability from below can be shown analogously. To this end we fix $t \in [0,T)$ and write 
 
 \[
  \begin{aligned}
   &\frac{I_F(t+h) -I_F(t)}{h}\\
   &= \frac{1}{h} \int_{t}^{t+h} F(t+h,s) \ \dd s + \int_{-\infty}^{t-2\varepsilon} \frac{F(t+h,s)- F(t,s)}{h} \ \dd s + \int_{t-2\varepsilon}^t \frac{F(t+h,s)- F(t,s)}{h} \ \dd s\\
   &=: \text{I}_{\rm{ii}}(h)+ \text{II}_{\rm{ii}}(h)+ \text{III}_{\rm{ii}}(h),
  \end{aligned}
 \]
where $\varepsilon>0$ is chosen according to assumption~\eqref{F_i_conditions}c). For the first term we derive by means of \eqref{F_i_conditions}a) that

\[
 \left|F(t,t)-\text{I}_{\rm{ii}}(h) \right| \le   \frac{1}{h} \int_{t}^{t+h}\left|F(t,t)- F(t+h,s) \right| \ \dd s \le \sup_{s \in [t,t+h]} \left|F(t,t)- F(t+h,s) \right| \to 0
\]

as $h \downarrow 0$. For the second term we use the mean value theorem for complex valued functions to see that

\[
 \begin{aligned}
 \int_{-\infty}^{t-2\varepsilon} \sup_{h \in (0, \varepsilon \land (T-t))} &\left|\frac{F(t+h,s)- F(t,s)}{h}\right| \ \dd s \le \int_{-\infty}^{t-2\varepsilon} \sup_{r \in [t, (t+\varepsilon) \land T]} \left| \dr F(r,s) \right| \ \dd s < \infty,
 \end{aligned}
\]

where the finiteness follows from~\eqref{e.diff_F_integral_bounded}. Therefore, we can apply the DCT to obtain

\[
 \lim_{h \downarrow 0} \text{II}_{\rm{ii}}(h) = \int_{-\infty}^{t-2 \varepsilon}  \dt F(t,s) \ \dd s.
\]

In order to tackle the third term we aim at applying the de la Vall\'ee-Poussin theorem again. For this purpose we let $1< q <(\beta+\gamma)^{-1}$ and deduce by using the mean value theorem for complex valued functions and~\eqref{e.diff_F_estimate} that

\[
 \begin{aligned}
   \sup_{h \in (0, \varepsilon \land (T-t))} &\int_{t-2\varepsilon}^t \left|\frac{F(t+h,s)- F(t,s)}{h} \right|^q \ \dd s \\
   &\le   \sup_{h \in (0, \varepsilon \land (T-t))}  \int_{t-2\varepsilon}^t \sup_{r \in [t,t+h]} \left|\dr F(r,s) \right|^q \ \dd s \\
   &\le \tilde C^q   \sup_{h \in (0, \varepsilon \land (T-t))}  \int_{t-2\varepsilon}^t \sup_{r \in [t,t+h]} |s|^{-\beta q}\left| r-s \right|^{-\gamma q}  \ \dd s \\
   &\le \tilde C^q   \int_{t-2\varepsilon}^t  |s|^{-\beta q}\left| t-s \right|^{-\gamma q}  \ \dd s \\
   &< \infty,
 \end{aligned}
\]

where the finiteness follows from the facts  $\beta q <1$ and $\gamma q <1$.
%
Now the de la Vall\'ee-Poussin theorem yields the uniform integrability of $\left( \frac{F(t+h,\cdot)-F(t,\cdot)}{h} \right)_{h \in (0, \varepsilon \land (T-t))}$ on $[t-2\varepsilon,t]$.  Consequently, we infer that

\[
 \lim_{h \downarrow 0} \text{III}_{\rm{ii}}(h) = \int_{t-2 \varepsilon}^t  \dt F(t,s) \ \dd s,
\]

which shows the differentiability from above with the desired right derivative. As mentioned above, the differentiability from below can be handled following the same line of argument.

By combining   \eqref{F_i_conditions_2} and \eqref{F_i_conditions} we conclude that every function $F$ fulfilling conditions \eqref{F_i_conditions}a)-\eqref{F_i_conditions}c) is continuously differentiable on $[0,\infty)$, since in this case \eqref{F_i_conditions_2} is applicable to the derivative $\dt F$.
\hfill$\square$

\section*{Acknowledgements}

Financial support by the Deutsche Forschungsgemeinschaft under grant BE3933/4-1 is gratefully acknowledged.



\end{document}